\documentclass[10pt]{article}

\usepackage{hyperref}
\usepackage{graphicx}
\usepackage{amssymb}
\usepackage{amsthm, amsmath}
\usepackage[all]{xy}
\usepackage{psfrag}

\def\g{\mathfrak{g}}        
\def\J{\mathbf{J}}          
\def\P{\mathcal{P}}     
\def\Rb{\mathbb{R}}     

\newtheorem{theorem}{Theorem}
\newtheorem{corollary}[theorem]{Corollary}
\newtheorem{definition}{Definition}
\newtheorem{proposition}{Proposition}
\newtheorem{lemma}{Lemma}
\newtheorem{rem}{Remark}
\newtheorem{example}{Example}

\begin{document}

\title{On the geometry of reduced cotangent bundles at zero momentum}

\author{Matthew Perlmutter\thanks{permanent address: Institute of Fundamental Sciences, Massey University, New Zealand.}, Miguel
Rodr\'{\i}guez-Olmos\thanks{corresponding author. permanent
address: Section de Mathematiques EPFL, CH-1015 Lausanne,
Switzerland. tel.: +41 216935507
fax: +41 216935839.},\\ M. Esmeralda Sousa-Dias\\
\small{Dep.
Matem\'atica, Instituto Superior T\'ecnico, Av. Rovisco Pais}\\
\small{1049--001 Lisboa, Portugal}}


\date{}

\maketitle

\begin{abstract}
We consider the problem of cotangent bundle reduction for proper
non-free group actions at zero momentum. We show that in this
context the symplectic stratification obtained by Sjamaar and
Lerman refines in two ways: (i) each symplectic stratum admits a
stratification which we call the secondary stratification with two
distinct types of pieces, one of which is open and dense and
symplectomorphic to a cotangent bundle; (ii) the reduced space at
zero momentum admits a finer stratification than the symplectic
one into pieces that are coisotropic in their respective
symplectic strata.
\end{abstract}

\section{Introduction}
This paper addresses the problem of symplectic reduction for
cotangent bundles with proper actions, at zero momentum. From the
point of view of mechanics, cotangent bundles are the most
important symplectic manifolds since they are the phase spaces for
most classical mechanical systems. The geometry of the reduced
space plays a crucial role in understanding the dynamics of
reduced Hamiltonian systems with non-freely acting symmetry
groups. We view this problem, then, as a fundamental one in the
theory of geometric mechanics and symplectic reduction.

A general theory of symplectic reduction for proper, and non-free
actions has been a subject of active research since the original
theory was worked out in Marsden and Weinstein \cite{MaWe} and
Meyer \cite{Me}. The geometric structure of the reduced spaces was
first satisfactorily understood, for the case of compact symmetry
groups, in the breakthrough paper of Sjamaar and Lerman
\cite{SaLe}, where the tools of stratification by orbit types were
first introduced to precisely determine how the reduced space,
which is not in general a manifold, is decomposed into symplectic
manifolds called symplectic strata. Indeed, from this point of
view, they were able to put in geometric context the earlier work
on this problem by proving that the symplectic strata of the
reduced space are the symplectic leaves of the reduced Poisson
algebra as determined in Arms et al. \cite{ArCuGo}. These
symplectic strata are obtained by first intersecting the zero
level set of the momentum map with the points in the original
symplectic manifold with the same orbit type, and then taking the
quotient of this space by the group action. They also explain how
the strata fit together by examining the behavior of a linear
symplectic action on a symplectic normal space, and applying the
Symplectic Slice Theorem due to Marle, Guillemin, and Sternberg.

Since this work, the field has continued to develop substantially.
In Bates and Lerman \cite{LeBa}, the theory was extended to proper
group actions and nonzero momentum, by way of orbit reduction,
with the assumption of locally closed coadjoint orbits. In Ortega
and Ratiu \cite{OrRa}, the theory of Poisson reduction by a free
Poisson action given in Marsden and Ratiu \cite{MaRa}, is extended
to the singular case. The symplectic reduction theory is extended
to the case of non-locally closed coadjoint orbits in Cushman and
\'Sniatycki \cite{CuSn} by looking at accessible sets of invariant
Hamiltonian vector fields. A comprehensive reference for all these
results, including several generalizations and improvements of the
theory and also their consequences in terms of reduction and
reconstruction of Hamiltonian dynamics is found in Ortega and
Ratiu \cite{OrRa2003}. Another text, Cushman and Bates
\cite{CuBa}, besides giving an overview of the general theory,
contains also many computed examples using invariant theory.

Specializing to cotangent bundles, one expects, as in the free
case, that the reduced space will admit special structure. Indeed,
in the free case, as is well known, the reduced space at zero
momentum is in fact simply the cotangent bundle of the orbit space
of the base with its canonical symplectic form. At nonzero
momentum it is known that the reduced symplectic space is
symplectomorphic to a coadjoint orbit bundle (see Marsden and
Perlmutter \cite{MarPer}). Alternatively it can be seen as the
image of a symplectic embedding into an appropriate cotangent
bundle (see for instance Marsden \cite{MarLec}).

Although various attempts were made to apply the general theory of
singular reduction to understand the important case of cotangent
bundles, until now, there has not been a complete picture without
strong assumptions. The literature begins with a result due to
Montgomery \cite{Mo} prior to the work of Sjamaar and Lerman in
which he extends the embedding theory of regular cotangent bundle
reduction to the case where the involved groups satisfy a special
dimension condition and the proper action on the base manifold is
assumed to consist of only one orbit type. In the paper
\cite{EmRo} Emmrich and R\"{o}mer give a complete solution to the
zero momentum reduced space for a proper action again with the
assumption that the base action consists of just one orbit type.
As one might guess from the free theory of cotangent bundle
reduction and the fact that the orbit space for the base action is
a manifold, they obtain that the reduced space at zero momentum is
just the cotangent bundle of the orbit space with its canonical
symplectic form.

The next paper to address the problem of reduction of cotangent
bundles is Lerman et al. \cite{LeMoSa}, where the example of $S^1$
acting on $T^{\ast}S^2$ is computed and the reduced space at zero
momentum is shown to be the ``canoe''. They also provide a result
for singular cotangent bundle reduction at zero in the case that
the action admits a cross section.

Finally, we note that in Schmah \cite{Schmah}, the results
obtained in \cite{EmRo} are again obtained with a different proof
and extended, under the same hypothesis on the isotropy groups in
the base, to deal with reduction at momentum values with trivial
coadjoint orbits.

\paragraph{The main results.}
There are several new results in this paper. An important guiding
principle in this work is that zero momentum reduced data should
correspond with data constructed from the group action on the
base, in particular the isotropy lattice.

We will consider a proper action of a Lie group $G$ on a manifold
$M$ and its lifted action on $T^*M$, which is Hamiltonian with
respect to the canonical symplectic form with equivariant momentum
map $\J:T^*M\rightarrow\g^*$. Our first main result, Theorem
\ref{teo3}, is that the isotropy lattice for the $G$-action on the
zero momentum level set, $\J^{-1}(0)$, is isomorphic to the
isotropy lattice for the base action of $G$ on $M$. We obtain
this, roughly, by decomposing $\J^{-1}(0)$ as a disjoint union of
fiber bundles along the base orbit types and then using a subtle
application of the Tube Theorem for slices. Next, relative to this
primary decomposition of $\J^{-1}(0)$, knowing its isotropy
lattice, we consider for each isotropy type $(L)$ the set
$(\J^{-1}(0))_{(L)}$. Call a pair of elements, $(H),(L)$ in the
isotropy lattice of $M$ a {\em connectable pair} over $(L)$,
provided $(H)\geq (L)$. This means that $L$ is conjugate to a
subgroup of $H$. Let us denote this relationship by $H\rightarrow
L$. We are now able to obtain a decomposition of the manifold,
$(\J^{-1}(0))_{(L)}$, into fiber bundles, one for each connectable
pair over $(L)$. That is, to each group larger than or equal to
$(L)$ in the base lattice we construct a fiber bundle, contained
in $(\J^{-1}(0))_{(L)}$.

The symplectic strata of the reduced space $\P_0:=\J^{-1}(0)/G$
are given by $\J^{-1}(0)_{(L)}/G$ for each $(L)$ in the base
isotropy lattice and we will further demonstrate that each of
these is in turn stratified by fiber bundles, which we call seams,
one for each connectable pair $H\rightarrow L$ over $(L)$. The
pair $L\rightarrow L$ is in fact identified under a natural
diffeomorphism to the cotangent bundle $T^{\ast}(M_{(L)}/G)$ and
we will prove that this is an open dense piece in this {\it
secondary} stratification of each symplectic stratum. The other
pieces fiber over the strata in the boundary of $M_{(L)}/G$.

The reduced symplectic structure fits together with respect to
this stratification in an elegant way. The cotangent bundle within
each symplectic stratum is open and dense. We prove in Proposition
\ref{secondprop} that the restriction of the reduced symplectic
form to each seam is in fact equal to the pull back of the
corresponding canonical symplectic form of the corresponding
cotangent bundle. In Theorem~\ref{secondsymp} we characterize the
reduced symplectic form on each symplectic stratum as the unique
extension of the canonical symplectic form of the open and dense
cotangent bundle (corresponding to the $L\rightarrow L$
connectable pair) to its closure. Furthermore, we prove in
Theorem~\ref{secondsymp} that the seams (corresponding to the
$H\rightarrow L$ pairs) are in fact coisotropic submanifolds
within their corresponding symplectic strata.

We consider the topology of the total reduced space $\P_0$ and
obtain a {\em co\-iso\-tro\-pic} stratification (Theorem
\ref{coisotropicdecomposition}) which demonstrates that the full
collection of objects, seams and cotangent bundles, corresponding
to the entire set of connectable pairs in the isotropy lattice of
the base, forms a stratification of $\P_0$, which is of second
order in the sense that each of its strata is labelled by a
connectable pair in the isotropy lattice. It is finer than the
stratification induced by the symplectic strata of Sjamaar and
Lerman and, in opposition to the latter, the continuous surjective
projection to $M/G$ happens to be a morphism of stratified spaces
with respect to the coisotropic stratification of $\P_0$ and the
orbit type stratification of $M/G$.

It should be noted that although the secondary and coisotropic
stratifications introduced here are suitable for explaining the
bundle structure of the reduced space for cotangent-lifted
actions, they lose some of the good properties enjoyed by the
symplectic stratification. First, the secondary and coisotropic
stratifications are not known to have the conical property, or to
satisfy the Whitney conditions as the symplectic stratification
(see \cite{SaLe, OrRa2003}). For general symplectic manifolds,
this is possible due to the existence of the Symplectic Slice
Theorem, which is not well adapted to the cotangent bundle case.
Second, unlike for the symplectic stratification, the secondary
pieces of the stratifications introduced here are not invariant
under the reduced Hamiltonian flows, which makes our results
difficult to apply to dynamics. However, these results do give
qualitative information about the evolution of isotropy in the
projection of the reduced flows to the reduced base space $M/G$.
See the end of the paper for more remarks on these two comments.

For most of the derivations of our results about these
stratifications we will work in the slightly weaker category of
$\Sigma$-decomposed spaces, because it is computationally simpler.
This category is introduced in Section 2. In Section~5, however,
we show how these results persist in the category of stratified
spaces.

\section{Background and preliminaries}

The main aim of this section is to review the results on proper
group actions and symplectic reduction that we shall need for the
rest of the paper. This review will also serve to fix notation. We
first review the basic results on proper group actions on
manifolds, namely the decomposition of the manifold into orbit
types which is a $\Sigma$-decomposition (to be introduced later)
of the manifold. We then recall the general theory of symplectic
reduction at zero momentum for proper group actions which
describes the decomposition of the reduced space at zero into
symplectic $\Sigma$-manifolds obtained in a natural way from the
orbit type decomposition of $\mathbf{J}^{-1}(0)$ (see \cite{SaLe}).
Finally, we will summarize the known results for cotangent bundle
reduction, first in the free case, and then, the next easiest case
for proper actions: the case with only one orbit type on the base
manifold.

\subsection{$\Sigma$--Decompositions and proper actions}

Recall that a smooth action of a Lie group on a manifold $M$ is
{\em proper} if the map $G\times M \rightarrow M\times M$,  $(g,
m) \mapsto (m, g\cdot m)$ is proper (the inverse image of a
compact set is compact).  Notice that we have denoted the action
map $G\times M\rightarrow M$ by a dot. For the proofs of the
following key properties see for instance  Duistermaat and Kolk
\cite{DuiKol} or Pflaum \cite{Pflaum}.

\bigbreak
\noindent{\it \bf Properties of proper actions:} Let $M$ be a
$G$-manifold with a proper action. Then,
\begin{enumerate}
\item The isotropy subgroup $G_m$ of any point $m\in M$ is
compact. \item Each orbit $G\cdot m$, $m\in M$, is a closed
submanifold of $M$ diffeomorphic to $G/G_m$. \item The orbit space
$M/G$ is  Hausdorff, locally compact  and paracompact. \item $M$
admits a $G$-invariant Riemannian metric. \item If all the
isotropy groups of points in $M$ are conjugate to a given one, the
orbit space $M/G$ is a smooth manifold and the projection
$M\rightarrow M/G$ is a surjective submersion.
\end{enumerate}

An important result for proper actions is the  standard model for
$G$-invariant neighborhoods. This is a consequence of the
existence of slices due to Koszul \cite{Koszul} in the case of $G$
compact and later extended to proper actions by Palais
\cite{Palais}. Let $\exp$ be  the exponential map associated to a
$G$-invariant metric and $S_m$ the orthogonal complement to
$\mathfrak{g}\cdot m=T_m(G\cdot m)$. Consider the product $G\times
S_m$ with the left diagonal action of $G_m$ given by $h\cdot
(g,v):=(gh^{-1},h\cdot v)$. This is well defined because by
construction $S_m$ is $G_m$-invariant. This action is free since
it is free in the first factor. Next, construct the associated
bundle $G\times_{G_m} S_m$ to the principal bundle $G\rightarrow
G/G_{m}$. There is a well defined $G$-action on this bundle given
by
$$
g\cdot[h,u]=[gh,u].$$
 With these constructions, one then has the
following result providing an explicit realization of a
$G$-invariant tubular neighborhood of the orbit through $m$.
\begin{theorem}[Tube Theorem]\label{slice}
The map $\phi:G\times_{G_m} S_m\rightarrow M$ given by
$$\phi([g,u])=g\cdot\mathrm{exp}_m(u)$$
restricts to a $G$-equivariant diffeomorphism from a $G$-invariant
neighborhood of the zero section of $G\times_{G_m} S_m$ to a
$G$-invariant neighborhood of $G\cdot m$ in $M$ satisfying
$$\phi([e,0])=m.$$
Consequently $\phi$ maps the set $[G,0]$, the
zero section of the bundle $G\times_{G_m}S_{m}$,
to the orbit $G\cdot m$.
\end{theorem}

\begin{rem}\label{remark1}{\rm
We can construct the $G$-invariant neighborhood of the zero
section of the associated bundle of the previous theorem as
follows. Let $r$ be some positive radius smaller than the
injectivity radius of $\exp_{m}$. Then the ball $B_r$ around $0$
in $S_{m}$ is $G_{m}$-invariant since the action is by isometries.
We refer to $B_r$ and $\exp_{m}(B_{r})$ as a linear slice and a
slice through $m$ for the $G$-action respectively. It is easy to
see that the $\exp_m$ map restricted to $B_r$ is a
$G_m$-equivariant diffeomorphism with respect to the linear action
of $G_m$ on $B_r$ and the base action of $G_m$ on $\exp_m(B_r)$
since the $G_m$ action must take geodesics to geodesics. Notice
then, that the only group elements leaving the slice invariant are
those in $G_m$, i.e. we have

\begin{equation}\label{propslice}
{\rm For\ any}\ z\in \mathrm{exp}_m\,(B_{r}),\ G_{z}\subseteq G_m.
\end{equation}
The $G$-invariant neighborhood of the zero section, alluded to in
the previous theorem, is then $G\times_{G_m}B_r$. The details of
the proof of the existence of slices for proper actions and of the
Tube Theorem can be found in \cite{DuiKol}.}$\hfill
\blacktriangle$
\end{rem}
\bigbreak
For a subgroup $H$ of a Lie group $G$ the conjugacy class of $H$
consists of all subgroups of $G$ that are conjugate to $H$ and will
be denoted by $(H)$. Denote by $I_M$ the set of conjugacy classes of
isotropy groups of points of $M$.  Corresponding to each element
of this set $(H)\in I_M$ we have the subset of $M$ of orbit type $(H)$
defined by
$$
M_{(H)}:=\{m\in M: G_{m}\in (H)\}.
$$
For a proper $G$ action on a manifold $M$ the connected components of the orbit type $M_{(H)}$ are embedded submanifolds.

In the set of conjugacy classes of $G$ we can define a partial
ordering $\leq$
  by $(H) \leq (K)$ if and only if $H$ is conjugate to a subgroup
of $K$ in $G$.  We will use the notation $(H)<(K)$ to mean that
$H$ is conjugate to a proper subgroup of $K$ in $G$, i.e. strictly
smaller than $K$. We will represent $I_M$ as a lattice in the
following way:  we draw an arrow  from  $H$ to $K$ when $H$ and
$K$ are representatives of two classes in $I_M$ such that
$(H)<(K)$ and there is no other class $(L)\in I_M$ such that
$(H)<(L)<(K)$.

For proper actions on a connected manifold $M$, Duistermaat and
Kolk~\cite{DuiKol}
 show the existence of a unique minimal class in the
isotropy lattice, say $(H_0)$. The orbit type
$M_{(H_0)}$ is called the principal orbit type and is open
and dense in $M$.

When a proper $G$-action on $M$ is not free then in general $M/G$
is not a manifold. It is usually said that $M/G$ is a stratified
space, with the strata being the sets $M_{(H)}/G$. It is so, of
crucial importance to our work to clarify the notion of
stratification by orbit types and most of our work we will done in
the weaker notion of a $\Sigma$--decomposition by the reasons
explained below. A comprehensive reference on the subject is
Pflaum \cite{Pflaum}.

Very often in the literature one encounters the stratification
notion as a decomposition of a topological space into pieces
(strata) that are manifolds satisfying the so-called frontier
condition (if $R\cap \overline{S}\neq\emptyset$ then $R\subset
\overline{S}$, for pieces $R,S$). As the following example from
Sjamaar and Lerman \cite{SaLe} shows, this stratification notion
is not adequate if we want to include $M/G$ as a stratified set
with strata $M_{(H)}/G$ since the set $M_{(H)}$, and consequently
$M_{(H)}/G$ is not in general a manifold, but a disconnected union
of manifolds of different dimensions.
\begin{example}
Consider the action of $S^1$ on $\mathbb{C}P^2$ given by
$$
e^{i\theta}\cdot [z_0,z_1,z_2]:=[e^{i\theta} z_0,z_1,z_2].
$$
It is clear that the orbit type submanifold $M_{(S^1)}$ is then the disjoint
union of the point at infinity $[1,0,0]$ and the complex plane $[0,z_1,z_2]$.
\end{example}

One could try to remedy this situation of the failure of $M_{(H)}$ to be a manifold by
considering a decomposition with pieces the connected components of
$M_{(H)}$. However in this case is not clear how the frontier conditions work. For these reasons we will
adopt here the notion of a $\Sigma$-decomposition.

\begin{definition}[$\Sigma$--decomposition]  Let $M$ be a paracompact
Hausdorff space with countable topology and ${\mathcal Z}$ a
locally finite partition of $M$ into locally closed subspaces
$S\subset M$. The pair $(M,{\mathcal Z})$ is called a
$\Sigma$-decomposed space and ${\mathcal Z}$ a
$\Sigma$-decomposition if the following conditions are satisfied:
\begin{itemize}
\item[i)] Every piece $S\in {\mathcal Z}$ is a $\Sigma$-manifold
in the induced topology, that is $S$ is a topological sum of
countably many connected smooth and separable manifolds. \item
[ii)] If $R\cap \overline{S} \neq \emptyset$, for a pair of pieces
$R, S \in {\mathcal Z}$, then $R\subset \overline{S}$ ({\it
frontier condition}).
\end{itemize}
\end{definition}
\paragraph{$\Sigma$-geometry.}
In general, a $\Sigma$-manifold will not be a manifold unless all
its connected components have the same dimension, however one can
reproduce virtually all the geometric results traditionally stated
for manifolds for these objects. In this sense, the tangent (resp.
cotangent) bundle $TM$ (resp. $T^*M$) of a $\Sigma$-manifold $M$
will be the topological sum of the tangent (resp. cotangent)
bundles of each connected component of $M$ and it is naturally a
$\Sigma$-manifold. A map $f:M\rightarrow N$ between
$\Sigma$-manifolds is smooth if the image of the intersection of
the domain of $f$ with each connected component of $M$ is
contained in a connected component of $N$ and the restriction of
$f$ to each connected component of $M$, seen as a map between
connected manifolds, is smooth. This allows us to implement the
concepts of diffeomorphisms, immersions, embeddings, etc of
$\Sigma$-manifolds. In the same spirit one can define vector
fields, flows, group actions, etc. Because of this flexibility,
many times we will simply drop the prefix $\Sigma$ when these
constructions arise, if the meaning is clear from the context.

The definition of a $\Sigma$-decomposition is well adapted to the decomposition of a $G$-manifold into orbit types.
Indeed, using   the Tube Theorem one can show that for a compact subgroup
$H$ of $G$ the sets $M_{(H)}$ are locally closed $\Sigma$-submanifolds of $M$,
meaning that each connected component of $M_{(H)}$ is a submanifold
of $M$ (for the proof see Corollary 4.2.8 and Lemma 4.2.9 of
Pflaum~\cite{Pflaum}). Furthermore  one can show that, for proper
actions,  the decomposition of $M$ into the $\Sigma$-submanifolds
$M_{(H)}$, is locally finite (see Pflaum~\cite{Pflaum} Lemma 4.3.2).
 We then have the following

 \begin{proposition}\label{prop1}
Let $M$ be a proper $G$-manifold.
  The orbit type decomposition of $M$
is  a $\Sigma$-decomposition with the pieces given by the orbit types
$M_{(H)}$, $(H)\in I_M$.
 In particular, the frontier condition for the pieces becomes equivalent to
 \begin{equation}\label{front}
 M_{(H)}\cap \overline{ M_{(K)}}\neq \emptyset \Longleftrightarrow
(K)\leq (H).
 \end{equation}
 \end{proposition}
Notice that the larger the orbit type, the smaller the isotropy
subgroup, that is  $(H)\leq(K)$ if and only if $M_{(K)}\subset \overline{M}_{(H)}$.

An useful way to visualize the global distribution of pieces of a
$\Sigma$-decomposed space $M$ is to associate to it a
decomposition lattice, where the elements are the pieces of $M$,
together with arrows showing the frontier conditions of pairs of
pieces. In this way, if $R$ and $S$ are two pieces we draw an
arrow from $R$ to $S$ if $R\subset \partial S$ and there is no
other piece $T$ such that $R\subset \partial T$, and $T\subset
\partial S$ where $\partial S:=\overline S \backslash S$. For
instance if our $\Sigma$-decomposition is the orbit type
decomposition of a $G$-manifold $M$, we find from the previous
proposition that the decomposition lattice of $M$ has the same
shape as the isotropy lattice of $I_M$, where in place of the
representative $H$ of an isotropy class we will have the
corresponding orbit type $M_{(H)}$, and the directions of the
arrows will be the reverse of those in the isotropy lattice.
Sometimes these particular kinds of decomposition lattices are
called orbit type lattices.

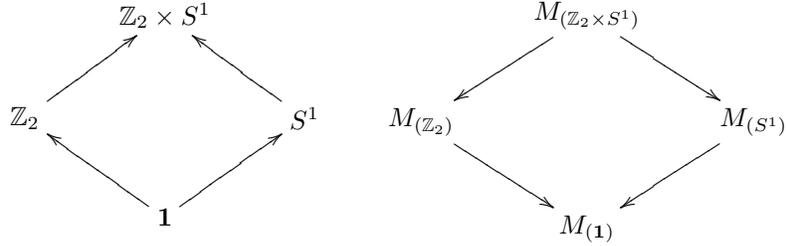
\begin{figure} \begin{center}
$$\xymatrix{ & \mathbb{Z}_2\times S^1  &  \\
\mathbb{Z}_2 \ar[ur] & & S^1 \ar[ul]\\
 & \mathbf{1} \ar[ul] \ar[ur] &
}\qquad
\xymatrix{ & M_{ \left(\mathbb{Z}_2\times S^1\right)} \ar[dl] \ar[dr] &  \\
M_{\left(\mathbb{Z}_2\right)} \ar[dr] & & M_{(S^1)} \ar[dl]\\
 &M_{(\mathbf{1})}&
}$$ \caption{\label{examples}Isotropy lattice and decomposition
lattice for the $\mathbb{Z}_2\times S^1$ action on
$M=\mathbb{R}^3$.}
    \end{center}\end{figure}

As an example consider the action of $\mathbb{Z}_2\times S^1$ on
$\Rb^3$ where $S^1$ acts by rotations around the $x_3$-axis and
$\mathbb{Z}_2$ by reflections with respect to the $(x_1,x_2)$-
plane. Since this group is compact, its resulting action on
$\mathbb{R}^3$ is proper and the isotropy groups are of four
types. $\mathbb{Z}_2\times S^1$ is the stabilizer of $\mathbf{0}$,
$\mathbb{Z}_2$ is the stabilizer of points of the $(x_1,x_2)$-
plane away from the origin, $S^1$ is the stabilizer of points of
the $x_3$-axis except the origin and the identity $\mathbf{1}$ is
the stabilizer of the remaining points. The respective isotropy
lattice and decomposition lattice are given in
Figure~\ref{examples}.

The $\Sigma$-decomposition of $M$ by orbit types   induces a
 $\Sigma$-decomposition on  $M/G$  (see for instance Theorem 4.3.10 of
Pflaum~\cite{Pflaum}). Its pieces  are $M_{(H)}/G$ where $H\in
I_M$ (recall that by item ({\it 5}) of the properties of proper
group actions these spaces are $\Sigma$-manifolds) and they
satisfy identical frontier conditions as the corresponding
$M_{(H)}$, so the decomposition lattices of $M$ and $M/G$ are
identical.

For further reference we define a morphism of decomposed spaces as follows.
\begin{definition}\label{defmorph}
A continuous map $f:(M_1, {\mathcal Z}_1)\rightarrow (M_2,
{\mathcal Z}_2)$ between decomposed spaces is called a morphism of
decomposed spaces if, for every piece $S\in {\mathcal Z}_1$ there
is a piece $R\in  {\mathcal Z}_2$ such that: $i)$ $f(S)\subset R$
and $ii)$ The restriction of $f$ to $S$ is smooth.

If all the restrictions $f\vert_S$ are injective, surjective,
immersions, submersions, embeddings etc, $f$ will be called   a
decomposed immersion, submersion, embedding, etc.

Finally, if $(M,\mathcal{Z}_1)$ and $(M,\mathcal{Z}_2)$ are two
decompositions of the same topological space $M$, we say that
$(M,\mathcal{Z}_1)$ is finer provided the identity map
$\mathrm{id}:(M,\mathcal{Z}_1)\rightarrow (M,\mathcal{Z}_2)$ is a
morphism of decomposed spaces.
\end{definition}
As a consequence of this definition, if $S_1$ and $S_2$ are two
pieces in $\mathcal{Z}_1$ whose images under $f$ are contained
respectively in $R_1$ and $R_2$ in $\mathcal{Z}_2$ and
$S_1\subset\overline{S_2}$ then $R_1\subset\overline{R_2}$.

\subsection{Symplectic reduction at zero momentum}
We now consider the setting of a Lie group $G$ acting properly and
symplectically on a symplectic manifold $\P$ and admitting an
equivariant momentum map $\J$. It has long been known since 1973,
1974 (in \cite{Me}, \cite{MaWe}) that when this action is free,
one can construct reduced symplectic manifolds
$\mathbf{J}^{-1}(\mu)/G_\mu$, henceforth referred to as
Marsden-Weinstein (MW) reduced spaces.

When the assumption of freeness of the action is removed, the
situation becomes immediately complicated as the momentum level
sets are no longer in general submanifolds. Nevertheless, with the
idea of partitioning the level sets into orbit types, it is
possible to prove that one can obtain a symplectic stratification
of the singular reduced spaces. In \cite{SaLe} the
Marsden-Weinstein reduced space at zero momentum    $\P_0 =
\J^{-1}(0)/G$, is described as a $\Sigma$-decomposition with each
piece a symplectic $\Sigma$-manifold constructed using orbit
types. In Theorem \ref{SL} we recall this result.

Throughout this paper we will use the following notations. Given a $G$-invariant subset $A$ of a $G$-manifold $\P$ we define
$$
A_{(H)} := A \cap \P_{(H)},\ \mathrm{and} \ \   A^{(H)} := A_{(H)}/G.
$$
We also make use of the following subsets of a $G$-manifold $M$:
$$M_H:=\{m\in M : G_m = H\},\qquad M^H:=\{m\in M : H\subset G_m\}.$$
Note that $M_{(H)}= G\cdot M_H$.

\begin{theorem}[Sjamaar and Lerman \cite{SaLe}] \label{SL}Let $(\P,
\omega)$ be a connected symplectic manifold on which $G$ acts
properly and symplectically admitting an equivariant  momentum map
$\J: \P\rightarrow \g^{\ast}$. Then
$\left(\J^{-1}(0)\right)_{(L)}$ is a $G$-invariant
$\Sigma$-submanifold of $\P$ and $\P_0 := \J^{-1}(0)/G$ is a
disjoint union of smooth symplectic $\Sigma$-manifolds,
\begin{equation}\label{sympstr}\P_0= \bigsqcup_{(L)\in I_\P} \P_0^{(L)},\end{equation} where
$\P_0^{(L)}:=\left(\J^{-1}(0)\right)_{(L)}/G$ with the reduced
symplectic form $\omega_0^{(L)}$ on $\P_0^{(L)}$ given by
$$\pi^{(L)*}\, \omega_0^{(L)}= i^{\ast}_{(L)}\,\omega,$$ where
$i_{(L)} : \left(\J^{-1}(0)\right)_{(L)}\rightarrow \P$ is the
inclusion, and the orbit projection is denoted by $\pi^{(L)}:
(\J^{-1}(0))_{(L)}\rightarrow \P_0^{(L)}$. Furthermore, this
partition of $\P_0$ is a $\Sigma$-decomposition with frontier
conditions obtained from the isotropy lattice $I_\P$.
\end{theorem}

\begin{rem}{\rm
In the above decomposition, some of the $\P_0^{(L)}$ might be
empty (this happens if $\J^{-1}(0)\cap \P_{(L)}=\emptyset$). We
will refer to \eqref{sympstr} as the {\em symplectic
decomposition} of $\P_0$.}$\hfill \blacktriangle$
\end{rem}

In the rest of the paper we study the additional structure that
the spaces $\P_0$ and $\P_0^{(L)}$ inherit from the cotangent
bundle structure of the original symplectic manifold $\P$
extending the known classical results for the free case.

\subsection{Cotangent bundle reduction}
In this section we review the well known results on cotangent
bundle reduction at zero momentum. We start with the free case and
then we review the case of a base manifold with just one orbit
type. Throughout this section we assume that $G$ is a Lie group
acting properly on a smooth manifold $M$ and by cotangent lifts on
$T^{\ast}M$.

The action of $G$ on $T^{\ast}M$ is
Hamiltonian with respect to the canonical symplectic form $\omega$
and has an $\operatorname{Ad}^{\ast}$-equivariant  momentum map $\J:
T^{\ast}M\rightarrow \g ^{\ast}$ given by
\begin{equation}
\label{mom}
\left\langle\J (p_m), \xi\right\rangle = \left\langle p_m,
\xi_M(m)\right\rangle,
\end{equation}
where $p_m\in T_{m}^{\ast}M$ and $\xi_M$ denotes the infinitesimal generator
for the $G$-action on $M$ corresponding to $\xi\in\g$.

In the free case, the cotangent lifted action on $T^{\ast}M$ is
also free and proper and consequently both orbit spaces, $M/G$ and
$T^{\ast}M/G$, are smooth manifolds. From \eqref{mom} one has
$$\alpha_{m}\in \mathbf{J}^{-1}(0)\cap T^{\ast}_{m}M \iff
\left\langle \alpha_{m},\xi_{M}(m)\right\rangle =0,
$$
and so the zero level set of $\J$ is   the annihilator of the bundle
$V\subset TM$ defined by $V_{m}=\{\xi_M(m) : \xi\in\g\} = T_m(G\cdot m)$.
That is, $\J^{-1}(0) = V^0$, which is a subbundle of $T^{\ast}M$.
The MW-reduced space $\P_0 := \J^{-1}(0)/G$ is
a smooth symplectic  manifold  with symplectic form $\omega_0$ induced
from the canonical symplectic form $\omega$ on $\P = T^*M$ defined by
$$\pi^*\omega_0 = i^{\ast}\omega,
$$
where $i:\J^{-1}(0)\rightarrow T^{\ast}M$ is the inclusion and
$\pi: \J^{-1}(0)\rightarrow \P_0$ the orbit projection map. The
following theorem,  due to Satzer in the case of $G$ Abelian, and
Abraham and Marsden in the general case shows that $\P_{0}$ is
symplectomorphic to the cotangent bundle of the orbit space $M/G$,
with its canonical symplectic form.
\begin{theorem}[Satzer \cite{Sat77}, Abraham and Marsden \cite{AbMa}]\label{teo1}
If $G$ acts freely  and properly on M and by cotangent lifts on
$T^*M$ then the symplectic reduced space $(\P_{0}, \omega_{0})$ is
symplectomorphic to $T^{\ast}(M/G)$ equipped with its canonical
symplectic structure.
\end{theorem}

\begin{proof}
We sketch a proof as follows. Consider the map
$\phi:TM/G\rightarrow T(M/G)$, defined by $\phi([v_{m}]) =
T_{m}\pi (v_{m})$. This  map is well defined and both fiber
preserving and surjective.  Its dual, $\phi^{\ast} : T^{\ast}(M/G)
\rightarrow T^{\ast}M/G$ is then a fiberwise injective bundle map
and Im$(\phi^{\ast}) = V^0/G$. As the vector bundles
$T^{\ast}(M/G)$ and $V^0/G$ over $M/G$ have the same dimension it
follows that $\phi^{\ast}$ is a bundle isomorphism, i.e
$T^{\ast}(M/G)\cong V^0/G$.  Finally, the symplectomorphism of the
theorem is given by  $(\phi^{\ast})^{-1}$.
\end{proof}
The next easiest generalization of this result, without the
freeness assumption, is the case where $M$ consists of a single
orbit type. This problem has been solved by Emmrich and R\"{o}mer
\cite{EmRo}, and later by Schmah \cite{Schmah} with a different
proof.

\begin{theorem}[Emmrich and R\"{o}mer \cite{EmRo}]\label{SOT} Let $G$ be a Lie
group acting on $M$ properly and on $\P = T^{\ast}M$ by cotangent
lifts. If all the points of $M$ have isotropy groups conjugate to
some $H \subset G$ (so that $M=M_{(H)}$), then
$\J^{-1}(0)=\left(\J^{-1}(0)\right)_{(H)}$  and $\P_0=(
\P_0)^{(H)}=\J^{-1}(0)/G$ is symplectomorphic to
$T^{\ast}(M/G)=T^{\ast}M^{(H)}$ with its canonical symplectic form
which we denote $\omega_{H}$.
\end{theorem}
\begin{rem}{\rm
 The symplectomorphism of the above theorem is the same as in
Theorem~\ref{teo1} for the free case.}$\hfill \blacktriangle$
 \end{rem}

 \section{Decomposition of $\J^{-1}(0)$}

In this section we prove a main result, Theorem~\ref{teo3},
showing that the isotropy lattice for the $G$-action on
$\J^{-1}(0)$ is identical to the isotropy lattice for the
$G$-action on the base manifold $M$. This result is special for
zero momentum and relies crucially on the fact that the zero
momentum level set corresponds to the annihilator of the tangent
spaces to the group orbits. Throughout the rest of the paper the
setting will be of a Lie group $G$ acting properly on a connected
smooth manifold $M$ and by cotangent lifts on $\P=T^*M$. Note that
the resulting action on $\P$ is automatically proper.

\subsection{Partition of $T^{\ast}M$ along orbit types}

 Due to the properness of the action, Proposition~\ref{prop1} gives that
 $M$ is a $\Sigma$-decomposed  manifold by orbit types,
that is
\begin{equation}\label{mot}
M = \bigsqcup_{(H)\in I_M}\, M_{(H)},
\end{equation}
where $M_{(H)}$ are $\Sigma$-submanifolds of $M$ verifying the frontier condition \eqref{front}.

Let $\mathbf{g}$ be a $G$-invariant
metric on $M$, and use \eqref{mot} to write $TM$
as a union of Whitney sums of $\Sigma$-vector bundles,  that is
\begin{equation}\label{tmot}
TM = \bigsqcup_{(H)\in I_M} TM_{(H)} \oplus NM_{(H)},
\end{equation}
where $NM_{(H)}$ denotes the orthogonal complement to $TM_{(H)}$ as a
$\Sigma$-vector bundle over $M_{(H)}$.

Since $G$ acts by isometries, t
he Legendre map
$\mathbb{F}L:TM\rightarrow T^{\ast}M$ defined by
$\mathbb{F}L(v_m)(w_m)=\mathbf{g}(m)(v_m,w_m),$ is an equivariant
bundle diffeomorphism from $TM$ to $T^{\ast}M$ and
  induces the
following dual splitting
\begin{equation}\label{tmp}
T^{\ast}M = \bigsqcup_{(H)\in I_M} T^{\ast}M_{(H)} \oplus N^{\ast}M_{(H)},
\end{equation}
which is a partition of $T^{\ast}M$.

Let $\J: T^{\ast}M\rightarrow\g^{\ast}$ be the $\operatorname{Ad}^{\ast}$-equivariant
momentum map for the
cotangent lifted action of $G$ on $T^{\ast}M$. The
partition \eqref{tmp} of  $T^{\ast}M$ along orbit types allows us to
express the zero level set of the momentum map as a disjoint union of
$\Sigma$-bundles over each orbit type in the base manifold.

\begin{proposition} \label{teo2}For a proper action of $G$ on the base
manifold $M$ the zero level set of the momentum map $\J$ for the
lifted $G$-action on $T^{\ast}M$ is a disjoint union of $\Sigma$-vector bundles
over $M_{(H)}$, where $(H)$ runs in the isotropy lattice $I_M$ of the
base manifold.  In particular
\begin{equation}\label{zls}
\J^{-1}(0) = \bigsqcup_{(H)\in I_M} \J^{-1}_{(H)} (0) \oplus N^{\ast}M_{(H)},
\end{equation}
where $\J_{(H)}$ is the momentum map for the $G$-action restricted  to the $\Sigma$-bundle $T^{\ast}M_{(H)}$ and $N^{\ast}M_{(H)}$ is
the $\Sigma$-conormal bundle of $M_{(H)}$.
\end{proposition}

\begin{proof}
Let $m \in M_{(H)}$ with stabilizer $G_m = H$. Recall that by definition of the momentum
map  \eqref{mom} we have
$$\J_m^{-1}(0) = \left(\g\cdot m\right)^\circ\subset T^{\ast}_mM,
$$
where we use the notation $\J_{m}:=\J\vert_{T_{m}^{\ast}M}$.
We will now decompose this annihilator making use of the metric $\mathbf{g}$
and the slice construction as follows. By definition of the normal bundle $NM_{(H)}$ to the
$\Sigma$-manifold $M_{(H)}$,  we have
\begin{equation}
\label{normal bundle}
T_mM=T_mM_{(H)}\oplus N_{m}M_{(H)}.
\end{equation}
Next, we use the metric to construct a linear slice $S_m$ for the action
of $G$ on $M$ at the point $m$,
$$
T_mM = \g \cdot m \oplus S_m,
$$
where $S_m$ is the orthogonal complement of the vertical space at $m$, i.e. $S_m = (\g\cdot m)^\perp = \left(T_m (G \cdot m)\right)^\perp$. We can decompose this space as follows noting
that $N_mM_{(H)}$ is orthogonal to $\g\cdot m\subset T_mM_{(H)}$,
$$
S_{m}=S_m\cap T_{m}M_{(H)}\oplus N_{m}M_{(H)}.
$$
Let us denote by $S_{m}':=S_{m}\cap T_{m}M_{(H)}$. Note that $S_{m}'$ is the
orthogonal complement in $T_mM_{(H)}$ to the subspace $\g\cdot m$. Therefore, by
construction, it is a linear slice for the $G$-action restricted to the $\Sigma$-manifold
$M_{(H)}$ through $m$. Consider the linear $H$ action on $S_{m}'$. Since $M_{(H)}$ has one orbit type by construction, $H$ must fix the entire space $S_{m}'$. In fact, letting $S_{m}^{H}$ denote the vector subspace of $S_m$ fixed by the $H$ action, we have $S_{m}^{H}=
S_{m}'$. To see this, if $(a,b)\in S_{m}^{H}\oplus N_{m}M_{(H)}$ is fixed by $H$
then $\exp_{m}\vert_{S_m}(a,b)\in M_H$ which implies that $(a,b)\in T_{m}M_{(H)}$ from which
we conclude that $b=0$.  We have therefore shown that $S_{m}^{H}=S_{m}'$,
and therefore we have the decompositions
\begin{equation}
\label{total}
T_{m}M=\g\cdot m\oplus S_{m}^{H}\oplus N_{m}M_{(H)}
\end{equation}
and
\begin{equation}
\label{partial}
T_{m}M_{(H)}=\g\cdot m\oplus S_{m}^{H}.
\end{equation}
Taking the dual of equation (\ref{total}) we obtain
$$
T^{\ast}_{m}M=(S_{m}^{H}\oplus N_mM_{(H)})^\circ\oplus (\g\cdot
m)^\circ \simeq (\g\cdot m)^{\ast}\oplus (S_{m}^H\oplus
N_mM_{(H)})^{\ast}
$$
so that $(\g\cdot m)^\circ\simeq (S_{m}^{H})^{\ast}\oplus
N^{\ast}_{m}M_{(H)}$. Furthermore, taking the dual of equation
(\ref{partial}) we obtain
\begin{equation*}
T^{\ast}_{m}M_{(H)}=(S_{m}^{H})^\circ\oplus \mathrm{Ann}(\g\cdot
m;T^{\ast}_mM_{(H)}) \simeq (\g\cdot m)^{\ast}\oplus
(S_{m}^{H})^{\ast}
\end{equation*}
so that
\begin{equation}
\label{annihilator}
\mathrm{Ann}(\g\cdot m;T^{\ast}_mM_{(H)})\simeq (S_{m}^{H})^{\ast}.
\end{equation}
Here we used the following notation: if $A\hookrightarrow B$ is a
linear injection of vector spaces, $\mathrm{Ann}(A;B^*)$ denotes
the annihilator of $A$ in $B^*$. Now, since the $G$-action
restricts to $M_{(H)}$ we can consider its cotangent lifted action
to $T^{\ast}M_{(H)}$. The momentum map for this action is just the
restriction of the momentum map on $T^{\ast}M$ to
$T^{\ast}M_{(H)}$. We call this momentum map
$\mathbf{J}_{(H)}:T^{\ast}M_{(H)}\rightarrow \g^{\ast}$. It then
follows from equation (\ref{annihilator}) that
$(S_{m}^{H})^{\ast}$ is the zero level set of the momentum map
$\J_{(H)}$ restricted to the fiber over $m\in M_{(H)}$. Denoting
by $\J_{(H)m}:=\J_{(H)}\vert_{T_{m}^{\ast}M_{(H)}}$ we have then
shown that
\begin{equation}\label{pointdecomposition}
\J_m^{-1}(0) = \J_{(H)m}^{-1}(0) \oplus N^{\ast}_mM_{(H)},
\end{equation}
from which the result follows.
\end{proof}

\subsection{Orbit types of $\J^{-1}(0)$}

In order to carry out the symplectic reduction for the zero level
set $\J^{-1}(0)$, Theorem \ref{SL} tells us that we need to
characterize $\P_0^{(L)}= \left(\J^{-1}(0)\right)_{(L)}/G$, for
each $(L)$ in the isotropy lattice of the $G$-lifted action on
$T^{\ast}M$.

By its very definition, the cotangent lifted action $G\times
T^{\ast}M\rightarrow T^{\ast}M$  satisfies $\tau (g\cdot
p_{m})=g\cdot \tau (p_{m})$ where
 the dot denotes both the left action of $G$ on $T^{\ast}M$ and on
$M$, and  $\tau :T^{\ast}M\rightarrow M$ denotes the projection.
It is then clear that in general the isotropy lattice for the
cotangent bundle, say $I_{T^*M}$, has more classes than $I_M$,
although it always contains those belonging to $I_M$ since $M$ is
$G$-equivariantly embedded in $T^{\ast}M$ as the zero section. The
main aim of this section is to show, in Theorem \ref{teo3}, that
there exists a one-to-one correspondence between orbit types in
$M$ and the symplectic pieces of the reduced space $\P_0 =
\J^{-1}(0)/G$. This is a remarkable feature of the zero momentum
level set. We start with the following coarse description of
$\J^{-1}(0)$ which will be refined in the subsequent theorem.

\begin{proposition} \label{prop2}The orbit types of the zero level set of the
momentum map for the cotangent lifted action of $G$ on $T^{\ast}M$ are
expressed as
\begin{equation}\label{eqprop2}
\left(\J^{-1}(0)\right)_{(L)} = \bigsqcup_{(H)\geq (L)} \J^{-1}_{(H)} (0)
\times \left(N^{\ast}M_{(H)}\right)_{(L)},
\end{equation}
where $(H)$ is in $I_M$ and $(L)$ is fixed in $ I_{T^{\ast}M}$.
\end{proposition}

\begin{proof}
As the projection $\tau:T^{\ast}M\rightarrow M$ is equivariant and
$(L)\in I_{T^{\ast}M}$, then $\tau((\mathbf{J}^{-1}(0))_{(L)})\cap
M_{(H)}\neq\emptyset$ implies $(L)\leq (H)$. So, from
\eqref{zls} we get
$$\left(\J^{-1}(0)\right)_{(L)} = \bigsqcup_{(H)\geq (L)}
\left(\J^{-1}_{(H)}
(0) \oplus N^{\ast}M_{(H)}\right)_{(L)}.
$$
Recall that $\J_{(H)}$ is the momentum map for the cotangent
lifted $G$-action to $T^{\ast}M_{(H)}$. We can now apply the
single orbit type theorem for cotangent lifted actions
(Theorem~\ref{SOT}) to obtain $\J^{-1}_{(H)} (0) =
\left(\J^{-1}_{(H)} (0)\right)_{(H)}$, which gives the result.
\end{proof}

At this point, we are able to get more information on the possible
subgroups $(L)$ by a careful analysis of the $G$-action on the
conormal bundles $N^{\ast}M_{(H)}$. The key to get finer
information is to apply the slice construction and the Tube
Theorem both for the $G$-action on $M_{(H)}$ and for the
$G$-action on $M$. This will allow us to relate the orbit types
for the $G$-action on the conormal bundle to the orbit types for
the $G$-action on the base. Specifically we find,
\begin{theorem}\label{teo3} For any $m \in M_{(H)}$ such that $G_m = H$, and any fixed $(L)\in I_{T^{\ast}M}$, then the orbit type $(L)$ of
the zero level set of the momentum map for the lifted $G$-action, restricted to the fiber over $m$, verifies
$$\left(\J^{-1}(0)\right)_{(L)}\cap T^{\ast}_mM \neq \emptyset$$
if and only if  both of the following  conditions hold
\begin{equation*}
i)\,  (L)\leq (H),\qquad\qquad ii)\,  M_{(L)}\neq\emptyset.
\end{equation*}

\end{theorem}

Before proving Theorem~\ref{teo3} we will prove a lemma relating the orbit types for the linear action of a
 subgroup $H$ of $G$ on $S_m = (\g\cdot m)^\perp$ and the orbit types of $G$ on the base manifold $M$.
It seems that most of the results in this lemma are scattered in the literature in a different form and so we present here a version that is better adapted to our purposes.

\begin{lemma} \label{lema1}Let $m \in M_{(H)}$ with $G_m = H$,
$B_r$ a ball of radius $r$ around
zero in $S_m = (\g\cdot m)^\perp$ with $r$ smaller than the injectivity
radius of $\exp_{m}$, $U$ and $\phi$ respectively the
$G$-invariant neighborhood of $G\cdot m$ and the diffeomorphism
given by the Tube Theorem and $M_{(K)}\neq \emptyset$ for some
$(K)\in I_M$. Consider the linear $H$ action on $S_m$. Then:
\begin{enumerate}
\item   $U \cap M_{(K)} \neq \emptyset $   if and only if
$(K)\leq(H)$.
\item
$
(S_m)_{(L)}\neq \emptyset,
$ if and only if there exists a class $(K)\in I_M$ with $(K)\leq (H)$ such
that $L$ is conjugate in $G$ to a representative  of $(K)$.
\item The set of points $[G,u]\subset G\times_{H}B_r$ with $u\in (B_r)_{(L)}$
gets mapped by $\phi$
into $M_{(K)}$ where $K$ is a subgroup of $H$ conjugate in $G$ to $L$.
\end{enumerate}
\end{lemma}

\begin{proof}  {\it 1.}:
Suppose $(K)  \leq (H)$, then by the frontier condition we have
$M_{(H)}\subset \overline{M_{(K)}}$. So, every open set in $M$
containing a point in $M_{(H)}$ must have nonempty intersection
with $M_{(K)}$.

Conversely, suppose  $m'\in U\cap M_{(K)}$. Then $G_{m'}$ is
conjugate to $K$ in $G$, i.e $(G_{m'})=(K)$. On the other hand,
$m'\in U$ and $U=G\cdot \exp_m(B_r)$, so $m' = g \cdot s$ for some $s\in
\exp_m(B_r)\subset M$ and $g\in G$.  Thus, as $m'=gs$ then $(G_{m'})=(G_s)$
and as $s\in\mathrm{exp}_m\,(B_r)$ then Remark~\ref{remark1} gives $G_s\subseteq H$. So $(G_{m'})=(G_s)=(K)\leq (H)$.

For {\it 2.}: From {\it 1,} we know that for  $(K)\leq (H)$ there
exists $s\in  \exp_m(B_r)$ such that $L:=G_{s}$ is conjugate to
$K$. Since $\exp_{m}$ is $H$-equivariant, the point
$\exp_{m}^{-1}(s)\cap S_m\in B_r$ is stabilized by $L$ under the
linear $H$ action on $B_r$. Since this action extends linearly to
the entire space $S_m$, we conclude that $(S_m)_{(L)}\neq
\emptyset$ where $L$ is conjugate to $K$. Conversely, let $L$ be a
subgroup of $H$ such that $(S_m)_{(L)}\neq \emptyset$. By
linearity of the $H$ action, $(B_r)_{(L)}\neq \emptyset$ and by
equivariance of $\exp_{m}$, we have $( \exp_m(B_r))_{(L)}\neq
\emptyset$. By {\it 1.}, this immediately implies that $L$ is
conjugate to $K$ for some $(K)\in I_M$ with $(K)\leq (H)$.

Finally, to prove {\it 3,}  it is sufficient to take $u\in (B_r)_{(L)}$ so that $H_u=L$. Now,
since $\exp_m$ is $H$ equivariant, we have that $\exp_m(u)$ is stabilized by $L$ as well and
in fact $G_{\exp_m(u)}=L$. It follows that each point in the set $\phi([G,u])=G\cdot \exp(u)$ is contained
in $M_{(L)}=M_{(K)}$ as required.

\end{proof}
\begin{proof} ({\it of Theorem~\ref{teo3}}) Recall from the proof of Proposition~\ref{teo2} that
$$\J_m^{-1}(0)= (\g\cdot m)^\circ \simeq S_m^{\ast} = (S_m^H \oplus N_mM_{(H)})^{\ast} \simeq (S_m^H)^{\ast}\oplus N_m^{\ast}M_{(H)}.$$

Since $H$ acts by isometries on $T_mM$ and on $T_mM_{(H)}$ by
restriction then $H$ maps $N_mM_{(H)}$ into itself and the action
of $H$ on $S_m^H\oplus N_mM_{(H)}$ is diagonal. Furthermore the
$H$ action on $S^H_m$ is trivial by construction.

Therefore for $(a,b)\in S_m^H \times N_mM_{(H)}$ one has
$H_{(a,b)} = H_b$, as $H_{(a,0)} = H$. Consequently, the orbit
type sets for the $H$ action on $S_m$ are of the form $S_m^H
\times \left(N_mM_{(H)}\right)_{(L)}$ where $(L)$ belongs to the
isotropy lattice for the linear $H$ action on $S_m$.

Let us  show that if $b\neq 0$ then $H_b$ is strictly contained in
$H$. For this,  note   that locally $S_m^H$ and $S_m$ are linear
slices  at $m$ for the $G$-actions on $M_{(H)}$ and $M$
respectively.

Consider  the direct product of $H$-invariant neighborhoods
$B_{r_1}\times B_{r_2}\subset S_m^{H}\oplus N_{m}M_{(H)}$ each of
them inside the disk of radius $r_i>0$ centered at $0$ in the
corresponding vector space, where $r_{1}^2+r_{2}^2<r^{2}$. Then,
their direct product is contained in the disk   $B_r\subset S_m$.
Denote by $\phi_{M_{(H)}}:G\times_{H}B_{r_1}\rightarrow
U_{M_{(H)}}$ the diffeomorphism from Theorem~\ref{slice} applied
to the slice for the $G$-action on $M_{(H)}$. The image of
$\phi_{M_{(H)}}$ is an open $G$-invariant set of $M_{(H)}$ and not
of $M$. Next consider the slice at the point $m$ for the entire
manifold $M$, modelled on the space $G\times_{H}(B_{r_1}\times
B_{r_2})$, and the corresponding map
$\phi:G\times_{H}(B_{r_1}\times B_{r_2})\rightarrow U$. Suppose
there exists $0\neq y\in B_{r_2}$ such that $H_y=H$. Then, the
entire open set $B_{r_1}\times ty$ where $t\in (0,r_2/||y||))$ is
stabilized by $H$ and therefore, by {\it 3}) of Lemma~\ref{lema1},
$\phi([G\times_{H}(B_{r_1}\times ty)])$ is contained in $M_{(H)}$.
However $\phi$ is a diffeomorphism so this image has one higher
dimension than $\phi_{M_{(H)}}(G\times_{H}B_{r_1})$. On the other
hand, they are both open sets in $M_{(H)}$, which is a
contradiction. We have then proved that $H_b\subsetneq H$ for
$b\neq 0$.

From {\it 2}) of Lemma~\ref{lema1} we know that $(S_m)_{(L)}\neq
\emptyset$ if and only if $L$ is conjugate to $K\subseteq H$ for
some $(K)\in I_M$ and $M_{(K)}\neq \emptyset$. Then we have proved
that
$$\left(\J_m^{-1}(0)\right)_{(L)} = S_m^H \times \left(N_m^{\ast}M_{(H)}\right)_{(L)} \neq \emptyset $$
if and only if
$$(L)\leq (H)\ \ \mbox{and}\ \  M_{(L)}\neq \emptyset.$$
\end{proof}

From the proof of Theorem~\ref{teo3} and noting that $M_{(H)} = G\cdot M_H$ we have
\begin{corollary}\label{cor1}
$\left(N^{\ast}M_{(H)}\right)_{(L)}\neq \emptyset$ if and only if
$(H)\geq (L)$ and $M_{(L)}\neq \emptyset$. Furthermore
$\left(N^{\ast}M_{(H)}\right)_{(H)}$ is the zero section of the
$\Sigma$-bundle $N^{\ast}M_{(H)} \rightarrow M_{(H)}$, i.e., it is
isomorphic to $M_{(H)}$.
\end{corollary}

To end this section we summarize in the next proposition the main results obtained so far for the orbit types of the zero momentum level set.

\begin{proposition} \label{prop10}In the previous conditions we have:
\begin{itemize}
\item[a)] $(L)\in I_{\J^{-1}(0)} \iff (L)\in I_M$ and then $
\mathcal{P}_0^{(L)}\neq \emptyset \iff (L)\in I_M$.
\item[b)] The
cotangent bundle projection $\tau$ restricts to the
$G$-equivariant continuous surjection $\tau_L :
\left(\J^{-1}(0)\right)_{(L)}\rightarrow \overline{M_{(L)}}$.

\item[c)] A fixed orbit type $(L)$ in the zero momentum level set
is a $\Sigma$-submanifold  of $T^*M$ which admits the following
$G$-invariant partition:
\begin{equation}\label{partition}
\left(\J^{-1}(0)\right)_{(L)} =\J^{-1}_{(L)}(0)\bigsqcup_{(H)>(L)} \J^{-1}_{(H)}(0) \times \left(N^{\ast}M_{(H)}\right)_{(L)}.
\end{equation}

\item[d)] For every $(H)>(L)$, the restrictions
$$t_L:={\tau_L}\vert_{\J^{-1}_{(L)}(0)} \quad\mbox{ and}
\quad t_{H\rightarrow L}:={\tau_L}\vert_{\J^{-1}_{(H)}(0)\times
\left(N^{\ast}M_{(H)}\right)_{(L)}} $$ are $G$-equivariant smooth
surjective submersions respectively onto $M_{(L)}$ and $M_{(H)}$.

\end{itemize}
\end{proposition}

\begin{proof}
Statement {\it a}) is proved in Theorem~\ref{teo3}. For {\it b}):
To prove continuity of $\tau_L$, first note that
$\overline{M_{(L)}}$ has the relative topology from $M$ so we must
show that for any open set $U$ in $M$, $\tau_{L}^{-1}(U\cap
\overline{M_{(L)}})$ is open in $(\J^{-1}(0))_{(L)}$. The
cotangent projection $\tau:T^{\ast}M\rightarrow M$ is of course
continuous, so $\tau^{-1}(U)$ is open in $T^{\ast}M$ and therefore
$\tau^{-1}(U)\cap (\J^{-1}(0))_{(L)}$ is an open set in
$(\J^{-1}(0))_{(L)}$. It is easy to show that, $\tau^{-1}(U)\cap
(\J^{-1}(0))_{(L)}=\tau_{L}^{-1}(U\cap \overline{M_{(L)}})$ from
which continuity of $\tau_{L}$ follows. $G$-equivariance is
obvious. Finally, note that the image of $\tau$ restricted to
$(\J^{-1}(0))_{(L)}$ is the disjoint union $\bigsqcup_{(H)\geq
(L)}M_{(H)}=\overline{M_{(L)}}$ since for each $(H)\geq(L)$,
$\J^{-1}_{(H)}(0) \times \left(N^{\ast}M_{(H)}\right)_{(L)}$ is a
$\Sigma$-fiber bundle over $M_{(H)}$. {\it c}) just follows from
Proposition \ref{prop2} and Theorem \ref{teo3}. To obtain {\it
d}), note that $\J_{(L)}^{-1}(0)$ is a $\Sigma$-fiber bundle over
$M_{(L)}$, i.e. disjoint union of smooth fiber bundles over each
connected component of $M_{(L)}$ and on each connected component
the fiber bundle projection $t_{L}$ is a smooth surjective
submersion. $G$-equivariance follows from the definition of the
cotangent lifted action. Similarly
$\J_{(H)}^{-1}(0)\times(N^{\ast}M_{(H)})_{(L)}$ is a
$\Sigma$-fiber bundle over $M_{(H)}$ with smooth surjective
$\Sigma$-submersion $t_{H\rightarrow L}$.
\end{proof}

\section{Topology and symplectic geometry of $\P_0$}
The general symplectic reduction theory (Theorem~\ref{SL}) tells
us that $\P_0$ is a $\Sigma$-decomposed space with symplectic
pieces $\P_0^{(L)}$. In the specific case of a cotangent bundle,
we show in this section that these symplectic pieces also admit a
$\Sigma$-decomposition which we call the {\em secondary
decomposition}. The pieces of the secondary decomposition of
$\P_0^{(L)}$ are studied in detail and we are able to prove that
there exists an open and dense piece which is diffeomorphic to the
cotangent bundle of $M_{(L)}/G$. The other pieces will be called
seams.

The reduced symplectic data then have a natural interpretation.
The reduced symplectic form $\omega_0^{(L)}$ in the symplectic
piece $\mathcal{P}_0^{(L)}$ can be obtained as the unique smooth
extension from this open dense part of the canonical symplectic
form on $T^{\ast}(M_{(L)}/G)$. Relative to the reduced symplectic
forms we will prove that the seams are coisotropic
$\Sigma$-submanifolds of $(\mathcal{P}_0^{(L)},\omega_0^{(L)})$.

We already know that the reduced space at zero momentum $\P_0$,
admits a symplectic $\Sigma$-decomposition in symplectic pieces
(Theorem \ref{SL}). We will prove that, joining together all the
pieces of the secondary decomposition of each symplectic piece
$\P_0^{(L)}$, the resulting partition of $\P_0$ is another
$\Sigma$-decomposition, which we call the {\em coisotropic}
decomposition. We explicitly identify the frontier conditions for
both $\Sigma$-decompositions of $\P_0$ and $\P_0^{(L)}$ and show
that the referred seams play a ``stitching role'', i.e. they
stitch the cotangent bundles appearing in the coisotropic
decomposition of $\P_0$, as we shall show in Theorem
\ref{coisotropicdecomposition}.

\subsection{The secondary decomposition of  $\P_0^{(L)}$}

We introduce the following notation. Recall that a connectable pair $H\rightarrow L$ is a pair of elements
$(H),(L)\in I_M$ such that $(H)\geq (L)$. Define the following fiber bundles
\begin{equation}
\label{preseam}
s_{H\rightarrow L}:=\J^{-1}_{(H)}(0) \times \left(N^{\ast}M_{(H)}\right)_{(L)}\rightarrow M_{(H)}.
\end{equation}
where the index $H\rightarrow L$ runs over the set of connectable pairs over a fixed isotropy class $(L)$.
As this is a $G$-invariant piece in the $G$-invariant partition \eqref{partition} of $(\J^{-1}(0))_{(L)}$ , we can
quotient by the $G$-action to obtain
\begin{equation}\label{seam}
S_{H\rightarrow L}:=\pi^{H\rightarrow L}(s_{H\rightarrow L})=
\frac{\J^{-1}_{(H)}(0) \times \left(N^{\ast}M_{(H)}\right)_{(L)}}{G}
\end{equation}
where $\pi^{H\rightarrow L}:=\pi^{(L)}|_{s_{H\rightarrow L}}$. We
shall then call $S_{H\rightarrow L}$, which is a fiber bundle over
$M_{(H)}/G$,  a {\it seam} from $H$ to $L$, and $s_{H\rightarrow
L}$, the fiber bundle over $M_{(H)}$, a {\it pre-seam}.

We then have the following partition of $\P_0^{(L)} =
\left(\J^{-1}(0)\right)_{(L)} /G$:
\begin{equation}\label{phpartition}
\P_0^{(L)} = \J^{-1}_{(L)}(0)/G\bigsqcup_{(H)>(L)} S_{H\rightarrow L}.
\end{equation}

Note that from Proposition~\ref{prop10}-$a)$ the conjugacy classes  $(L)$ and $(H)$
appearing in the above equations belong to $I_M$, with $(L)$ fixed in the disjoint union.
Moreover, due to the $G$-equivariance of the restrictions of the cotangent bundle projection, referred to in $b)$ and $d)$ of Proposition~\ref{prop10}, we have
 \begin{itemize}\label{pageitem}
 \item[i)] The map $\tau_L$ descends to a continuous surjection, say $\tau^L : \P_0^{(L)}\rightarrow \overline{M^{(L)}}$,
 where $\overline{M^{(L)}}$ is the closure of $M^{(L)}=M_{(L)}/G$.
 \item[ii)] For every $(H)>(L)$, the maps $t_L$ and $t_{H\rightarrow L}$ of Proposition~\ref{prop10}-$d)$ descend to
the following surjective submersions
 $$t^L: \J^{-1}_{(L)}(0)/G\rightarrow M^{(L)}\qquad t^{H\rightarrow L}: S_{H\rightarrow L}\rightarrow M^{(H)}.
 $$
 \end{itemize}
 These maps are summarized in the following commutative diagrams.
 $$\xymatrix{ \mathbf{J}_{(L)}^{-1}(0)/G \ar@{^{(}->}[r]^{~i_0^{L}} \ar[d]_{t^{L}}   &
 \mathcal{P}_0^{(L)} \ar[d]^{\tau^{L}}  \\
             M^{(L)}   \ar@{^{(}->}[r]^{i^{L}}         &
             \overline{M^{(L)}}}\quad\mathrm{and}\quad
 \xymatrix{ S_{H\rightarrow L} \ar@{^{(}->}[r]^{~i_0^{H\rightarrow L}} \ar[d]_{t^{H\rightarrow L}}
 & \mathcal{P}_0^{(L)} \ar[d]^{\tau^{L}}  \\
             M^{(H)}   \ar@{^{(}->}[r]^{i^{H}}         &
             \overline{M^{(L)}}}$$
 Note that   we know, from the general symplectic reduction theory, that $\P_0^{(L)}$ is a
 smooth (symplectic) $\Sigma$-manifold, but,  recalling that $M^{(L)}:=M_{(L)}/G$,   $\overline{M^{(L)}}$ in general is only a topological
 space, endowed with the relative topology of $M/G$.
 In the next proposition we show that $\overline{M^{(L)}}$ is a $\Sigma$-decomposed space
 and we identify the frontier conditions for its pieces.
 \begin{proposition} \label{prop4}
 $\overline{M^{(L)}}$ is a $\Sigma$-decomposed space with pieces $M^{(H)}$,
 for all $(H) \geq (L)$. The frontier conditions are given by
 $$M^{(K)}\cap \overline{M^{(H)}} \neq \emptyset \Longleftrightarrow (K)
 \geq (H).
 $$
 Furthermore $M^{(L)}$ is open and dense in $\overline{M^{(L)}}$.

 \end{proposition}
 \begin{proof}
 Using that
 $\overline{M_{(L)}}= \bigsqcup_{(H)\geq (L)} M_{(H)}$ and
 $$ \overline{M^{(L)}}= \bigsqcup_{(H)\geq (L)} M_{(H)}/G = \bigsqcup_{(H)\geq (L)} M^{(H)} .$$
 Since the orbit type decomposition of $M$ is a $\Sigma$-decomposition with pieces $M_{(H)}$, for all $(H)\in I_M$,
it is easy to see that $\overline{M_{(L)}}$ is also a $\Sigma$-decomposed space with
pieces $M_{(H)}$ with $(H)\in I_M$ and $(H)\geq (L)$. Since an orbit type decomposition of $M$ induces a
 $\Sigma$-decomposition of $M/G$ with pieces $M_{(H)}/G$ then, by the same argument as before,
 $\overline{M^{(L)}}$ is a $\Sigma$-decomposed space with the obvious frontier conditions stated in the proposition.

Therefore it remains to  prove that $M^{(L)}$ is open and dense in
$\overline{M^{(L)}}$. Density is obvious. For openness, consider a
point $x\in M^{(L)}=M_{(L)}/G$ and an open neighborhood $U'$ of
$x$ in $\overline{M^{(L)}}$. This means that there exists an open
neighborhood $U$ of $x$ in $M/G$ with $U'=U\cap
\overline{M^{(L)}}$. Adjusting $U$ we can assure that $U\cap
M^{(H)}=\emptyset$ for every $(H)>(L)$, since the points that are
stabilized by $(H)$ lie in the boundary of $M_{(L)}$. For such a
$U$ then, $U'=U\cap \overline{M^{(L)}}$ is totally contained in
$M^{(L)}$.
\end{proof}

The element $\J^{-1}_{(L)}(0)/G$ of the  partition~\eqref{phpartition} of $\P_0^{(L)}$ is
diffeomorphic to the cotangent bundle of $M^{(L)}$ by the single
orbit type theorem (Theorem~\ref{SOT}), since $\J_{(L)}$ is the
momentum map for the restriction of the $G$-action to
$T^{\ast}M_{(L)}$. We will denote this piece by $C_L$ and the
partition~\eqref{phpartition} can be written as
 \begin{equation}\label{phpartition2}
\P_0^{(L)} = C_L\bigsqcup_{(H)>(L)} S_{H\rightarrow
L},
\end{equation}
for all $(L),(H)\in I_M$. Note also that the piece $C_L$ of the
partition~\eqref{phpartition2}, which is diffeomorphic to a
cotangent bundle, can also be seen as a seam from $L$ to $L$
since, by Corollary~\ref{cor1},
$\left(N^{\ast}M_{(L)}\right)_{(L)}$ is the zero section of the
$\Sigma$-bundle $N^{\ast}M_{(L)} \rightarrow M_{(L)}$ and so
Definition \eqref{seam} gives
$$C_L = S_{L\rightarrow L} \simeq \J^{-1}_{(L)}(0)/G \simeq T^{\ast}M^{(L)}.
$$
If there is no danger of confusion we will use  $S_{L\rightarrow L}$,  $C_L$  and  $T^{\ast}M^{(L)}$ to denote the same piece.
Before stating the main result of this subsection we need to prove
the openness of the surjective map $\tau_L$ given in
Proposition~\ref{prop10}-$b) $
\begin{lemma}\label{lem2}
The  map $\tau_L: (\mathbf{J}^{-1}(0))_{(L)}\rightarrow \overline{M_{(L)}}$ is an open map.
In addition, the quotient map, $\tau^{L}:\P_{0}^{(L)}\rightarrow \overline{M_{(L)}}/G$ is
also open.
\end{lemma}
\begin{proof}
We begin by considering, for a fixed $(H)\geq (L)$ ,
$s_{H\rightarrow
L}:=\mathbf{J}_{(H)}^{-1}(0)\times(N^{\ast}M_{(H)})_{(L)}
\hookrightarrow T^{\ast}M |_{M_{(H)}}\hookrightarrow T^{\ast}M$.
The above sequence is then a sequence of embedded
$\Sigma$-submanifolds. Furthermore, the pre-seam $s_{H\rightarrow
L}$ is a $\Sigma$-fiber bundle over $M_{(H)}$ which embeds as a
$\Sigma$-fiber subbundle of the $\Sigma$-vector bundle
$T^{\ast}M|_{M_{(H)}}$. Since the topology of $(\J^{-1}(0))_{(L)}$
and $s_{H\rightarrow L}$ for each $(H)\geq (L)$ is the relative
topology of a $\Sigma$-submanifold of $T^{\ast}M$, the open sets
of  $(\J^{-1}(0))_{(L)}$ are $(\J^{-1}(0))_{(L)}\cap U$ for each
open set $U$ in $T^{\ast}M$. To prove the openness of the map
$\tau_L$ we need to show that $\tau_L(\J^{-1}(0)_{(L)}\cap U)$ is
an open set in $\overline{M_{(L)}}$. Now, since
\begin{equation}
\label{first step}\begin{array}{ccc}
\tau_L((\J^{-1}(0))_{(L)}\cap U) & = & \tau_{L}\left(\bigsqcup_{(H)\geq (L)}s_{H\rightarrow L}\cap U\right)\\
& = &
\bigsqcup_{(H)\geq (L)}t_{H\rightarrow L}(s_{H\rightarrow L}\cap U),
\end{array}\end{equation}
we need to consider the sets $t_{H\rightarrow L}(s_{H\rightarrow L}\cap U)$ contained in $M_{(H)}$.
In fact we will establish the following intersection formula for an arbitrary open set $U\subset T^{\ast}M$,
\begin{equation}
\label{intersection formula}
t_{H\rightarrow L}(U\cap s_{H\rightarrow L})=\tau(U)\cap M_{(H)},
\end{equation}
from which the proof  of openness will be an easy consequence.
Abstracting slightly, given an embedding of fiber bundles, where the
embeddings are inclusions,
$$\xymatrix{ A_1\ar@{^{(}->}[r] \ar[d]_{\tau_1}
 & A_2 \ar[d]^{\tau_2}  \\
             M_1   \ar@{^{(}->}[r]         &
             M_2}
$$
and given an open set $U$ in $A_2$, it is a general result that
$$\tau_2(U)\cap M_1=\tau_1(U\cap A_1).$$
Notice that since the fiber projection maps $\tau_1$ and $\tau_2$
are surjective submersions, they are open maps and therefore
the left hand side of the previous equation is open in $M_1$ since
its open sets are generated from the relative topology and $\tau_2(U)$ is
an open set in $M_2$. Similarly the right hand side is also an open
set in $M_1$. Note that this result also holds for a $\Sigma$-fiber bundle
embedding. Applying this result to the $\Sigma$-fiber bundle $s_{H\rightarrow L}\hookrightarrow T^{\ast}M$ which fibers over the base inclusion $M_{(H)}\hookrightarrow M$, we conclude
that the intersection formula (equation \eqref{intersection formula}) holds and therefore,
following equation \eqref{first step} we have,
\begin{align*}
\tau_L((\J^{-1}(0))_{(L)}\cap U)&=\bigsqcup_{(H)\geq (L)}t_{H\rightarrow L}(s_{H\rightarrow L}\cap U)=\bigsqcup_{(H)\geq(L)}\tau(U)\cap M_{(H)}\\
&=\tau(U)\cap \bigsqcup_{(H)\geq (L)}M_{(H)}=\tau(U)\cap \overline{M_{(L)}}.
\end{align*}
However, $\tau(U)\cap \overline{M_{(L)}}$ is an open set in
$\overline{M_{(L)}}$ since $\tau(U)$ is open in $M$ and
$\overline{M_{(L)}}$ has the relative topology from $M$.

Next we consider the map $\tau^{L}$ defined through the $G$-equivariance of the
map $\tau_{L}$ giving the following commutative diagram.
$$\xymatrix{ (\J^{-1}(0))_{(L)}\ar@{->}[r]^{\tau_{L}} \ar[d]_{\pi^{(L)}}
 & \overline{M_{(L)}} \ar[d]^{\pi_{M}^{(L)}}  \\
             \mathcal{P}_{0}^{(L)}  \ar@{->}[r]_{\tau^{L}}         &
             \overline{M^{(L)}}}
$$
The vertical arrows in this diagram are open maps since they are quotients
of a $G$-action and the topology on the base is given by the quotient topology.
Therefore, by openness of the map $\tau_L$, given an open set $U$ in $\P_{0}^{(L)}$, the set $\tau_{L}((\pi^{(L)})^{-1}(U))$ is
open in $\overline{M_{(L)}}$, and therefore since
\begin{equation*}
\tau^{L}(U)=\pi_{M}^{(L)}(\tau_{L}((\pi^{(L)})^{-1}(U)))
\end{equation*}
we conclude, by openness of the map $\pi_{M}^{(L)}$, that
$\tau^{L}(U)$ is open.
\end{proof}

We are now able to prove one of the main results of this section.

\begin{theorem}\label{thm8} The partition~\eqref{phpartition2} is a $\Sigma$-decomposition of $\P_0^{(L)}$ that will be
called the secondary decomposition of $\mathcal{P}_0^{(L)}$. The
piece $C_L$ is open and dense and diffeomorphic to
$T^{\ast}M^{(L)}= T^{\ast}(M_{(L)}/G)$. The frontier
conditions are:
\begin{itemize}
\item[1)] $S_{H\rightarrow L} \subset \partial\, C_L$ for all $(H) >(L)$.
\item[2)] $S_{H'\rightarrow L} \subset \partial S_{H\rightarrow L}$ if and only if $(H') >(H)>(L)$.
\end{itemize}
The map $\tau^L : \P_0^{(L)}\rightarrow \overline{{M}^{(L)}}$ is a
$\Sigma$-decomposed surjective submersion.
\end{theorem}
\begin{proof} By construction of \eqref{phpartition2}  and because an orbit type decomposition is a
    $\Sigma$-decomposition it is then clear that the partition~\eqref{phpartition2}
    is a locally finite partition. Since the pieces of the
partition are $\Sigma$-submanifolds of $\P_{0}^{(L)}$ then they are
automatically locally closed.

Let us prove that $C_L$ is open and dense. Let $U$ be an open
neighborhood of $z\in \P_0^{(L)}=\left(
\J^{-1}(0)\right)_{(L)}/G$. Since, by Lemma~\ref{lem2}, the map
$\tau^L: \P_0^{(L)}\rightarrow \overline{M^{(L)}}$ is open,
$\tau^L(U) = O$ is an open set in $\overline{M^{(L)}}$. By
Proposition~\ref{prop4}, $M^{(L)}$ is dense in
$\overline{M^{(L)}}$ and so $O\cap M^{(L)}\neq \emptyset$. For
$y\in O\cap M^{(L)}$, we have $(\tau^{L})^{-1}(y)=(t^L)^{-1}(y)
\subset C_{L}$ and $(\tau^L)^{-1}(y) \cap U \neq \emptyset$.  It
follows that, $U\cap C_L\neq \emptyset$, proving the density. For
the openness of $C_L$ note that by Proposition~\ref{prop4},
$M^{(L)}$ is open and so
$\left(\tau^L\right)^{-1}\left(M^{(L)}\right) =C_L$ is also open
by the continuity of $\tau^L$.

For {\it 1)}, let $z\in S_{H\rightarrow L} $ with $(L) <(H)$ and $U$ an open neighborhood
of $z$. As $C_L$ is dense in $\P_0^{(L)}$ then $U\cap C_L \neq \emptyset$.
Furthermore as $C_L$ and   $S_{H\rightarrow L} $ are disjoint for $(L) <(H)$ it follows that
$z\in \partial C_L$.

Let us now prove {\it 2)}.
By the openness property of $\tau^L$ then any neighborhood $U$ of a point $z\in  S_{H' \rightarrow L}$ in $\mathcal{P}_0^{(L)}$ is
mapped by $\tau^L$ to an open
neighborhood of $\tau^{L}(z)$ in $\overline{M^{(L)}}$, say $O$.
Then $O\cap M^{(H)}\neq \emptyset$ if and only if $(H')>(H)>(L)$
because $\overline{M^{(L)}}$ is a $\Sigma$-decomposed space. Then for
$y\in O\cap M^{(H)}$ we have $(t^{H\rightarrow L})^{-1}(y)\cap U
\neq\emptyset$, proving {\it 2)}.

The map $\tau^L$ restricted to each seam is a surjective
submersion, that is $\tau^L(S_{H'\rightarrow L})= t^{H'\rightarrow
L}(S_{H'\rightarrow L}) = M^{(H')}$, also  $\tau^L(S_{H
\rightarrow L})= t^{H\rightarrow L}(S_{H\rightarrow L}) =
M^{(H)}$. By the frontier conditions we get that $\tau^L$ is a
$\Sigma$-decomposed surjective submersion.
\end{proof}

We  will now describe   the symplectic structure of
the symplectic pieces $\mathcal{P}_0^{(L)}$. Recall that by the
single orbit type theorem (Theorem~\ref{SOT}), for each $(H)\in I_M$  there is a
diffeomorphism
\begin{equation}\label{psih}
\psi^H:C_H\rightarrow T^{\ast}M^{(H)}
\end{equation}
which is a $\Sigma$-bundle map covering the identity in
$M^{(H)}$. Consider now,  for each piece in the
partition~\eqref{partition} of $(\mathbf{J}^{-1}(0))_{(L)}$  the
projection,
$$\mathrm{p}_{1\,H\rightarrow L}:\mathbf{J}^{-1}_{(H)}(0)\times
(N^{\ast}M_{(H)})_{(L)}\rightarrow \mathbf{J}^{-1}_{(H)}(0).$$
  Notice that this map is
  just the identity map on the first element of the partition,
$\mathbf{J}^{-1}_{(L)}(0)$. These are equivariant maps
that descend to surjective submersions
\begin{equation}\label{p1hl}
\mathrm{p}_1^{H\rightarrow L}:S_{H\rightarrow L}\rightarrow
C_H=\mathbf{J}_{(H)}^{-1}(0)/G.
\end{equation}
Then for any connectable pair $H\rightarrow L$ over $(L)$,  we have for the corresponding
piece $S_{H\rightarrow L}$ of $\mathcal{P}_0^{(L)}$, a surjective submersion
\begin{equation}\label{psibarhl}
\overline{\psi}^{\,H\rightarrow
L}=\psi^H\circ\mathrm{p}_1^{H\rightarrow L}:S_{H\rightarrow
L}\rightarrow T^{\ast}M^{(H)}
\end{equation}
covering the identity on $M^{(H)}$. In the particular case
$(H)=(L)$ we have that $S_{L\rightarrow L}=C_L$ and
$\overline{\psi}^{\,L\rightarrow L}=\psi^L$ is a diffeomorphism.
If we denote by $\omega_{H}$ the canonical symplectic form in
$T^{\ast}M^{(H)}$ we can then induce on each piece of the
secondary decomposition of $\mathcal{P}_0^{(L)}$ a closed two form
by
$$\begin{array}{l}
\mbox{on}\quad C_L:\,\, \Omega_L:={\psi^L}^{\ast}
\omega_{L},\quad\mbox{and}\quad \mbox{on}
\quad
S_{H\rightarrow L}: \,\, \Lambda_{H\rightarrow
L}:={\overline{\psi}^{\,H\rightarrow L}}^*\omega_{H}.
\end{array}$$
Then $\Omega_L$ is symplectic and $\Lambda_{H\rightarrow L}$ is
 degenerate.

By Theorem \ref{SL} the piece $\mathcal{P}_0^{(L)}$ has an
abstractly defined reduced symplectic form $\omega_0^{(L)}$. Is
then natural to ask to what extent the structures introduced so
far are compatible. The answer to this question is given in the
next proposition, which  together with Theorem~\ref{secondsymp}
are the main results characterizing the symplectic geometry of
$\mathcal{P}_0^{(L)}$.

\begin{proposition}\label{secondprop}
Let  $T^{\ast}M^{(H)}$ be equipped with the canonical symplectic
form $\omega_H$ and  $\P_0^{(L)}$ with the symplectic form $
\omega_0^{(L)}$ given by  Theorem~\ref{SL}. Let also
$\overline{\psi}^{\,H\rightarrow L}$ and $\psi^L$ be respectively
the surjective submersion \eqref{psibarhl} and the diffeomorphism
\eqref{psih}. Then, there are closed two forms $\Omega_L$ on $C_L$
and $\Lambda_{H\rightarrow L}$ on $S_{H\rightarrow L}$ defined by
$$\Omega_L:={\psi^L}^{\ast} \omega_L,\qquad \Lambda_{H\rightarrow
L}:=\left({\overline{\psi}^{\,H\rightarrow L}}\right)^*\omega_H,
$$
verifying
$$
{\it i)}\,\, \omega_0^{(L)}\vert_{C_L}=\Omega_L\qquad\mbox{and}
\qquad {\it ii)}\,\,\omega_0^{(L)}\vert_{S_{H\rightarrow
L}}=\Lambda_{H\rightarrow
L}.
$$
\end{proposition}
\begin{proof}
We will present the proof for {\it ii)} from which {\it i)}
follows by taking $(H)=(L)$ and noting that
${\overline{\psi}^{\,H\rightarrow L}} = \psi^L$ . First note that
by Theorem~\ref{SL} the symplectic form, $\omega_0^{(L)}$, in
$\P_0^{(L)}$ is given by
\begin{equation}\label{symp1}
\pi^{(L)^*} \omega_0^{(L)} =  i_{(L)}^*\omega,
\end{equation}
   where $\omega$ is
the canonical symplectic form in $T^*M$,  $\pi^{(L)}$ and
$i_{(L)}$ respectively the orbit projection and the inclusion
defined in the referred theorem (see also diagram below). In order
to prove equation {\it ii)} let us consider the following diagram
$$\xymatrix{
&s_{H\rightarrow L}\ar@{^{(}->}[r]^(0.4){i_{H\rightarrow L}}
\ar[d]^{\pi^{H\rightarrow L}} &\left(\J^{-1}(0)\right)_{(L)}
\ar@{^{(}->}[r]^(0.6){i_{(L)}}\ar@{->}[d]_{\pi^{(L)}}&T^*M\\
T^*M^{(H)}&S_{H\rightarrow L}
\ar@{->}[l]_(0.3){\overline{\psi}^{\,H\rightarrow L}}
\ar@{^{(}->}[r]^{i_0^{H\rightarrow L}}&\P_0^{(L)}&
}
$$
As $\pi^{H\rightarrow L}$ is a submersion, if we prove
\begin{equation}\label{sympext}
\left(\pi^{H\rightarrow L}\right)^*\, \left(i_0^{H\rightarrow
L}\right)^*\, \omega_0^{(L)}=
\left(\pi^{H\rightarrow L}\right)^*\, \,\left(
{\overline{\psi}^{\,H\rightarrow
L}}\right)^*\omega_{H},
\end{equation}
the claim  $ \left(i_0^{H\rightarrow L}\right)^*\, \omega_0^{(L)}=
\left( {\overline{\psi}^{\,H\rightarrow L}}\right)^*\omega_{H}$ of
the proposition  follows.

  From the above diagram we have
$i_0^{H\rightarrow L} \circ \pi^{H\rightarrow L} = \pi^{(L)} \circ
i_{H\rightarrow L}
$. So the left hand side of \eqref{sympext} becomes
\begin{equation}\label{symp2}
\left(\pi^{H\rightarrow L}\right)^*\, \left(i_0^{H\rightarrow
L}\right)^*\, \omega_0^{(L)}= i_{H\rightarrow L}^*\, \pi^{(L)*}\,
\omega_0^{(L)} =  i_{H\rightarrow L}^*\, i_{(L)}^*\omega,
\end{equation}
where the second identity follows from the definition  \eqref{symp1} of
$\omega_0^{(L)}$.

Note that the image of $i_{(L)} \circ i_{H\rightarrow L}$ is
contained in
$T^*M\vert_{M_{(H)}}\subset T^{\ast}M$. Therefore,  denoting by $\phi$ and
$i_H$ the following inclusions
$$\xymatrix{s_{H\rightarrow L} \ar@{^{(}->}[r]^(0.4){\phi}
&T^*M\vert_{M_{(H)}} \ar@{^{(}->}[r]^(0.6){i_H} &T^*M,}
$$
equation \eqref{symp2} is equivalent to
\begin{equation}\label{symp3}
\left(\pi^{H\rightarrow L}\right)^*\, \left(i_0^{H\rightarrow
L}\right)^*\, \omega_0^{(L)} = i_{H\rightarrow L}^*\,
i_{(L)}^*\omega = \phi^*\, i^*_H\, \omega.
\end{equation}
So in order to prove \eqref{sympext} it remains to show that
\begin{equation}\label{sympext2}
\phi^*\, i^*_H\, \omega = \left(\pi^{H\rightarrow L}\right)^* \,\left(
{\overline{\psi}^{\,H\rightarrow
L}}\right)^*\omega_H.
\end{equation}
For \eqref{sympext2} recall that $\overline{\psi}^{\,H\rightarrow
L}=\psi^H\circ\mathrm{p}_1^{H\rightarrow L}$. Then the right hand
side of  \eqref{sympext2} is given by
\begin{equation}\label{symp4}
\begin{array}{ll}
\left(\pi^{H\rightarrow L}\right)^* \,\left(
{\overline{\psi}^{\,H\rightarrow
L}}\right)^*\omega_H &= \left(\pi^{H\rightarrow L}\right)^*\,
\left(\mathrm{p}_1^{H\rightarrow L}\right)^*\, \left(\psi^H\right)^*
\omega_H\\
&\\
&=\left(\mathrm{p}_{1_{ H\rightarrow
L}}\right)^*\,\left(\pi^H\right)^*\, \left(\psi^H\right)^*\, \omega_H
\end{array}
\end{equation}
where the second identity follows from the commutativity of the
following diagram
$$\xymatrix{
s_{H\rightarrow L}\ar[d]_{\pi^{H\rightarrow L}} \ar[r]^{\mathrm{p}_{1_{
H\rightarrow L}}} &\J^{-1}_{(H)}(0)\ar[d]^{\pi^H} \\
S_{H\rightarrow L}\ar[r]^(0.4){\mathrm{p}_1^{H\rightarrow L}}
&\J^{-1}_{(H)}(0)/G
}
$$
Recall that $\J_{(H)}$ is the momentum map for the $G$-action on
$T^{\ast}M_{(H)}$ and so by the single orbit type theorem
(Theorem~\ref{SOT}), $\J^{-1}_{(H)}(0)/G$ is symplectic with symplectic form, $(\psi^{H})^{\ast}\omega_{H}$, induced from the canonical symplectic form
$\omega_{(H)}$ on $T^{\ast}M_{(H)}$, given by
\begin{equation}\label{can2}
(\pi^{H})^{\ast}\ (\psi^{H})^{\ast} \omega_{H} = j^{\ast} \omega_{(H)},
\end{equation}
where $j$ denotes the inclusion $j : \J^{-1}_{(H)}(0)\rightarrow
T^{\ast}M_{(H)}$.

Using equation \eqref{can2} and substituting into \eqref{symp4},
we obtain
\begin{equation}\label{symp5}
\left(\pi^{H\rightarrow L}\right)^{\ast} \,\left(
{\overline{\psi}^{\,H\rightarrow L}}\right)^{\ast}\omega_H =
\left(\mathrm{p}_{1_{ H\rightarrow L}}\right)^{\ast}\, j^{\ast}
\omega_{(H)}.
\end{equation}
The map $j\circ \mathrm{p}_{1_{ H\rightarrow L}}$ is related with
$\phi$ by $j\circ \mathrm{p}_{1_{ H\rightarrow L}} = p\circ \phi$ where
$p$ is the projection $p: T^*M\vert_{M_{(H)}} = T^*M_{(H)} \oplus
N^*M_{(H)}\rightarrow T^*M_{(H)} $. That is, we have the following
commutative diagram.
$$\xymatrix{
s_{H\rightarrow L}\ar@{^{(}->}[r]^(0.4){\phi}  \ar[d]_{\mathrm{p}_{1_{
H\rightarrow L}}} &T_{M_{(H)}}^*M\ar[d]^{p} \\
\J^{-1}_{(H)}(0)\ar@{^{(}->}[r]^(0.5){j} &T^*M_{(H)}}
$$
Therefore equation \eqref{symp5} is equivalent to
\begin{equation}\label{symp6}
\left(\pi^{H\rightarrow L}\right)^{\ast} \,\left(
{\overline{\psi}^{\,H\rightarrow
L}}\right)^{\ast}\omega_H = \left(\mathrm{p}_{1_{ H\rightarrow
L}}\right)^{\ast}\, j^{\ast}\,\omega_{(H)} = \phi^{\ast} \, p^{\ast}\,
\omega_{(H)}.
\end{equation}
So in order to finish the proof of \eqref{sympext2} it is sufficient to
show that
\begin{equation}\label{symp7}
p^{\ast} \, \omega_{(H)} = i^{\ast}_H \, \omega,
\end{equation}
which will be done  in local coordinates.

Let $({\mathcal U}, x_{1}, \cdots, x_{n})$ be a coordinate system
on $M$ adapted to $M_{(H)}$, so that ${\mathcal U}\cap M_{(H)}$ is
described by $x_{k+1}=\cdots =x_{n}=0$. Let $(T^{\ast}{\mathcal
U}, x_{1}, \cdots, x_{n}, \xi_{1}, \cdots , \xi_{n})$ be the
associated cotangent coordinate system on $T^{\ast}M$. Let
$\Theta$ and $\Theta_{(H)}$ be the canonical one-forms
respectively on $T^{\ast}M$ and on $T^{\ast}M_{(H)}$. In these
local coordinates the maps $i_{H}$ and $p$ are
$$i_{H}(x, \xi)= (x_{1}, \cdots , x_{k}, 0\cdots, 0, \xi_{1}, \cdots ,
\xi_{n})$$
$$p(x, \xi) = p(x_{1}, \cdots , x_{k}, 0\cdots, 0, \xi_{1}, \cdots ,
\xi_{n})= (x_{1}, \cdots , x_{k}, \xi_{1}, \cdots, \xi_{k})$$
Then,
$$p^{\ast}\Theta_{(H)}= \sum_{i=1}^k \xi_{i} dx_{i}
= \sum_{i=1}^{n}\xi_{i}\, dx_{i}\vert_{{\rm
span}\left\{\frac{\partial}{\partial x_{i}};\, 1\leq i\leq
k\right\}}= i_{H}^*\Theta ,$$
and the result \eqref{symp7} follows for the respective symplectic
forms by taking the exterior derivative.
\end{proof}

The previous proposition describes in part the abstract reduced
symplectic
form $\omega_0^{(L)}$ by means of natural explicitly constructed
closed two-forms on each piece of the secondary decomposition. However
this is not a complete description since we cannot say what is
$\omega_0^{(L)}$ at a point of a seam applied to vectors that are
not tangent to that seam. The next theorem gives a
characterization of the reduced form, as well as information on
the symplectic data of the $\Sigma$-submanifolds that form the
secondary
decomposition.

\begin{theorem}\label{secondsymp}
In the conditions of Proposition \ref{secondprop}, the reduced symplectic form
$\omega_0^{(L)}$ of the symplectic piece $\mathcal{P}_0^{(L)}$ is
the unique smooth extension of $\Omega_L$ from $C_L$ to
$\mathcal{P}_0^{(L)}$.
Furthermore, the following are satisfied:
\begin{enumerate}
\item $C_L$ is an open dense maximal symplectic
$\Sigma$-submanifold of $(\mathcal{P}_0^{(L)},\omega_0^{(L)})$
symplectomorphic to $(T^*M^{(L)},\omega_L)$
\item $S_{H\rightarrow L}$ are coisotropic $\Sigma$-submanifolds of $(\mathcal{P}_0^{(L)},\omega_0^{(L)})$
\end{enumerate}
\end{theorem}
\begin{proof}
Consider a point $x\in\mathcal{P}_0^{(L)}$ and two vectors $X_x,Y_x\in T_x\mathcal{P}_0^{(L)}$. Because $C_L$ is open and dense we can find
a sequence of points $x_k\in C_L$ and vectors $X_{x_k},Y_{x_k}\in T_{x_k}C_L\simeq T_{x_k}\mathcal{P}_0^{(L)}$ such that
$$
\mathrm{lim}\, x_k  =  x,\quad
\mathrm{lim}\, X_{x_k}  =  X_x, \quad
\mathrm{lim}\, Y_{x_k}  =  Y_x.
$$
We can then study the existence of the limit of the sequence $\Omega_L(x_k)(X_{x_k},Y_{x_k})$ as $k\rightarrow \infty$. By Proposition \ref{secondprop} we have that
$$\mathrm{lim}\,\Omega_L(x_k)(X_{x_k},Y_{x_k})=\mathrm{lim}\, \omega_0^{(L)}(x_k)(X_{x_k},Y_{x_k})=\omega_0^{(L)}(x)(X_x,Y_x)$$
where in the first equality we have used openness and density of
$C_L$ through the identification $T_{x_k}C_L\simeq
T_{x_k}\mathcal{P}_0^{(L)}$, and the last equality comes from
continuity of $\omega_0^{(L)}$. So, we have proved that there
exists a unique continuous extension of $\Omega_L$ to
$\mathcal{P}_0^{(L)}$. That this extension is smooth follows from
the fact that $\omega_0^{(L)}$ is the extension and is known to be
smooth by general reduction theory. The restrictions of this
extension to $C_L$ and to each seam follow tautologically from
Proposition \ref{secondprop}.

{\it 1}) is a trivial consequence of Theorem \ref{thm8} and Proposition \ref{secondprop}.
To prove {\it 2}), first recall from symplectic linear algebra (see \cite{LiMarle} for instance)
that for $(V,\omega)$ a
symplectic vector space and $W$ a vector subspace, then $W$ is
coisotropic if and only if $\mathrm{rank}\,(\omega\vert_W)=2\dim\,
W-\dim\, V$.

In our case we will do this dimension counting with respect to the
following tangent spaces. First fix $x \in S_{H\rightarrow
L}\subset \mathcal{P}_0^{(L)}$ and let $y\in s_{H\rightarrow L}$
be such that $x=\pi^{H\rightarrow L}(y)$ and $G_y=L$. Note that we
can always find such a $y$. Now denote by $z:=t_{H\rightarrow
L}(y)$ the projection of $y$ to the base $\Sigma$-manifold
$M_{(H)}$ so that $G_{z}\in (H)$. Let us call $H':=G_{z}$. Then we
set $V=T_x\mathcal{P}_0^{(L)},\, W=T_xS_{H\rightarrow L}$ and
$\omega=\omega_0^{(L)}(x)$. Note that by Proposition
\ref{secondprop} we have $\omega\vert_W=\Lambda_{H\rightarrow
L}(x)$.

Now, since $T^{\ast}M^{(L)}$ is open and dense in $\mathcal{P}_0^{(L)}$,
\begin{equation*}
\mathrm{dim}\, V=\mathrm{dim}\,T^{\ast}M^{(L)}=2\,(\mathrm{dim}\,
M_{(L)}-\mathrm{dim}\, G+\mathrm{dim}\, L).
\end{equation*}
On the other hand, by
construction of $\Lambda_{H\rightarrow L}$, we have that
\begin{equation*}
\mathrm{rank}\, \omega\vert_W=\mathrm{dim}\,
T^{\ast}M^{(H)}=2\,(\mathrm{dim}\, M_{(H)}-\mathrm{dim}\,
G+\mathrm{dim}\, H).
\end{equation*}
Finally, we have to compute $\mathrm{dim}\,
W=\mathrm{dim}\, S_{H\rightarrow L}$. For this, note that
$\mathrm{dim}\, S_{H\rightarrow L}=\mathrm{dim}\,(\mathbf{J}^{-1}(0)\cap
T_{z}^{\ast}M)_{(L)}+\mathrm{dim}\,
M_{(H)}-\mathrm{dim}\,G+\mathrm{dim}\,L$. Where $(L)$ refers to
the linear $H'$-action. On the other hand, the Legendre transform
maps $(\mathbf{J}^{-1}(0)\cap T^{\ast}_zM)_{(L)}$ $H'$-equivariantly
isomorphically to $(S_z)_{(L)}$. Now, if $\phi$ and $U$ are the
diffeomorphism and neighborhood of $z$ in $M$ given by the Tube
Theorem, then $\phi$ restricts to a diffeomorphism between
$G\times_{H'} (S_z)_{(L)}$ and $U\cap M_{(L)}$.
Since $\mathrm{dim}\,G\times_{H'}(S_{z})_{(L)}=\mathrm{dim}\,G+\mathrm{dim}\,(S_{z})_{(L)}-\mathrm{dim}\, H$,
we can compute
\begin{equation*}
\mathrm{dim}\,(S_z)_{(L)}=\mathrm{dim}\,
M_{(L)}-\mathrm{dim}\, G+\mathrm{dim}\, H.
\end{equation*}
Finally we obtain
$\mathrm{dim}\,W=\mathrm{dim}\,M_{(H)}+\mathrm{dim}\,M_{(L)}-2\mathrm{dim}\,G
+\mathrm{dim}\,H +\mathrm{dim}\,L$. It is then clear that the
condition $\mathrm{rank}\,(\omega\vert_W)=2\dim\, W-\dim\, V$ is
always satisfied.
\end{proof}
As a straightforward application of dimension counting we obtain
the following result
\begin{corollary}
We have the following facts about seams,
\begin{itemize}
\item[i)] If $(H)\neq (L)$ the seam $S_{H\rightarrow L}$ can never
be a symplectic submanifold of $\mathcal{P}_0^{(L)}$. \item[ii)]
If $(H)\neq (L)$, the seam $S_{H\rightarrow L}$ is a coisotropic
submanifold whose symplectic leaf space associated to the null
foliation of $\Lambda_{H\rightarrow L}$ is symplectomorphic to
$T^{\ast}M^{(H)}$ with its canonical symplectic form. \item[iii)]
A connected component of $S_{H\rightarrow L}$ is a Lagrangian
submanifold of $\mathcal{P}_0^{(L)}$ if and only if the
corresponding connected component (i.e. under the projection
$t^{H\rightarrow L}$) of $M^{(H)}$ is zero-dimensional.
\end{itemize}
\end{corollary}
\begin{proof} For {\it i}), It is obvious that if $H\neq L$ then
$\Lambda_{H\rightarrow L}$ has nonzero kernel. To see {\it ii}),
we know from the Theorem~\ref{secondsymp} that $S_{H\rightarrow
L}$ is a coisotropic submanifold of $\P_{0}^{(L)}$ and that the
restriction of the symplectic form to the seam satisfies
\begin{equation*}
\omega_0^{(L)}\vert_{S_{H\rightarrow
L}}=\left({\overline{\psi}^{\,H\rightarrow L}}\right)^*\omega_H,
\end{equation*}
where, recall, $\overline{\psi}^{\,H\rightarrow L}:S_{H\rightarrow L}\rightarrow T^{\ast}M^{(H)}$
is a surjective submersion. Since the symplectic leaf space is characterized by precisely
this equation, {\it ii}) follows.
For {\it iii}),
note that $S_{H\rightarrow L}$ is coisotropic, so it is
Lagrangian if and only if it has minimal dimension, i.e. $\frac
12\mathrm{dim}\,\mathcal{P}_0^{(L)}=\mathrm{dim}\,M^{(L)}=\mathrm{dim}\,M_{(L)}-
\mathrm{dim}\,G+\mathrm{dim}\,L$. Recalling from the proof of the
last theorem that $\mathrm{dim}\,S_{H\rightarrow
L}=\mathrm{dim}\,M_{(H)}+\mathrm{dim}\,M_{(L)}-2\,\mathrm{dim}\,G
+\mathrm{dim}\,H +\mathrm{dim}\,L$ we obtain that the Lagrangian
condition is satisfied if
$\mathrm{dim}\,M_{(H)}-\mathrm{dim}\,G+\mathrm{dim}\,H=0$, but this
is nothing but the dimension of $\mathrm{dim}\,M^{(H)}$.
\end{proof}

\subsection{The coisotropic decomposition of $\P_0$}
In this section we analyze the global structure of the topological
space $\mathcal{P}_0$, describing a new, cotangent-bundle adapted,
decomposition that is finer than the symplectic one. Recall from
previous sections that for each isotropy class $(L)$ in $M$ there
is a symplectic piece $\mathcal{P}_0^{(L)}$ in the reduced space
and the converse is also true. Furthermore, each of these pieces
is again a $\Sigma$-decomposed space with an open and dense piece
$C_L$, diffeomorphic to the cotangent bundle $T^{\ast}M^{(L)}$ and
a collection of seams $S_{H\rightarrow L}$, one for each
connectable pair $H\rightarrow L$ over $(L)$ satisfying $(H)\neq
(L)$. In this sense we obtained that the $(L)$-type symplectic
piece of the zero momentum reduced space has the structure of a
``topological fiber bundle'' over $\overline{M^{(L)}}$, where the
continuous projection $\tau^L$ is a $\Sigma$-decomposed surjective
submersion.

We want now to extend this bundle picture to the whole reduced
symplectic space $\mathcal{P}_0$. First of all, let
$\tau_0=\tau\vert_{\mathbf{J}^{-1}(0)}$ be the restriction of the
cotangent bundle projection to the zero momentum level set, which
is $G$-equivariant, and $\tau^0$ the corresponding descended map
$\tau^0:\mathcal{P}_0\rightarrow M/G$. By similar arguments to
those in the previous section, $\tau^0$ is a continuous surjective
open map. It should be immediately noticed that it is not a
morphism of $\Sigma$-decomposed spaces if $\mathcal{P}_0$ is
endowed with the symplectic decomposition and $M/G$ with the orbit
type one, since by Theorem \ref{thm8} the image of
$\mathcal{P}_0^{(L)}$ is contained in the closure of $M^{(L)}$ and
it has nonempty intersection with the boundary. It is our aim to
explain how a different decomposition of $\mathcal{P}_0$ in terms
of cotangent bundles and seams can be given in a way such that
$\tau^0$ is a $\Sigma$-decomposed surjective submersion. Consider
the following partition of $\mathcal{P}_0$:
\begin{equation}\label{coisotropicpartition}
\mathcal{P}_0=\bigsqcup_{(L)}C_L\bigsqcup_{(K')>(K)}S_{K'\rightarrow
K} \quad \mathrm{for\, all}\quad (L),(K),(K')\in I_M
\end{equation}
Obviously $\tau^0$ restricts on each piece to the previously defined
smooth surjective submersions
$$\tau^0\vert_{C_L}=t^L:C_L\rightarrow M^{(L)}\quad\mathrm{and}\quad
\tau^0\vert_{S_{K'\rightarrow K}}=t^{K'\rightarrow
K}:S_{K'\rightarrow K}\rightarrow M^{(K')}$$ The next theorem
explains the properties of this partition as well as the bundle
structure of $\mathcal{P}_0$.

\begin{theorem}\label{coisotropicdecomposition}
The partition \eqref{coisotropicpartition} of $\mathcal{P}_0$ is a
$\Sigma$-decomposition, that we will call coisotropic decomposition, and
satisfies:
\begin{enumerate}
\item If $(H_0)$ is the principal orbit type in $M$ then $C_{H_0}$ is
open
and dense in $\mathcal{P}_0$.
\item The frontier conditions are:
\begin{itemize}
\item[(i)] $C_{K}\subset \partial C_{H}\quad \mathrm{if\, and\, only\,
if}\quad
(H)<(K)$.
\item[(ii)] $S_{K\rightarrow H}\subset\partial C_{H}\quad \mathrm{if\,
and\,
only\, if}\quad  (H)<(K)$.
\item [(iii)] $C_{K}\subset\partial S_{K\rightarrow H}\quad
\mathrm{if\, and\,only\, if}\quad  (H)<(K)$.
\item [(iv)]$S_{K'\rightarrow H}\subset\partial S_{K\rightarrow
H}\quad\mathrm{if\, and\, only\, if}\quad  (H)<(K)<(K')$.
\item [(v)]$S_{K\rightarrow H'}\subset\partial S_{K\rightarrow
H}\quad\mathrm{if\, and\, only\, if}\quad  (H)<(H')<(K)$.
\end{itemize}
\item The continuous projection $\tau^0:\mathcal{P}_0\rightarrow
M/G$ is a $\Sigma$-decomposed surjective submersion with respect
to the coisotropic decomposition of $\mathcal{P}_0$ and  the usual
orbit type decomposition of $M/G$. \item If $I_M$ has more than
one class the coisotropic decomposition is strictly finer than the
symplectic decomposition, otherwise they are  identical.
\end{enumerate}
\end{theorem}

\begin{proof}
\noindent ({\it 1}) Note that by Proposition \ref{prop10}  the
symplectic $\Sigma$-decomposition of $\P_0$ has pieces
$\P_0^{(L)}$ for every $(L)\in I_M$. So, if $(H_0)$ is the
principal orbit type in $M$ then $\P_0^{(H_0)}$ is an open and
dense $\Sigma$-submanifold of $\P_0$. As $C_{H_0}$ is open and
dense in $\P_0^{(H_0)}$ with respect to the topology in
$\P_0^{(H_0)}$ and this topology is the  induced one from the
topology in $\P_{0}$, the result follows.

({\it 2}) The items {\it (ii)} and {\it (iv)} follow from Theorem
\ref{thm8} regarding the pieces in the statement as pieces of the
decomposition  of $\P_0^{(H)}$ for the respective $(H)$. Also,
{\it (iii)} follows from {\it (v)} by taking the limit $(H') =
(K)$.  Then, it remains to show {\it (i)} and {\it (v)}.

For {\it (i)}, recall that from the symplectic $\Sigma$-decomposition of $\P_0$ we
have the following frontier conditions
$$\P_0^{(K)} \subset \partial \P_0^{(H)} \Longleftrightarrow (H) < (K).
$$
    As $C_K\subset \P_0^{(K)}\subset  \partial \P_0^{(H)}$ if and only if
$ (H) < (K)$, then any open set $V_{x}$ in $\P_0$ containing a point
$x\in C_K$
must have nonempty intersection with $\P_0^{(H)}$ if and only if
$(H) < (K)$. But, since $C_H$ is dense in $\P_0^{(H)}$ it follows that
 $V_{x}$ also intersects $C_H$, proving  {\it (i)}.

For {\it (v)}: First note that   a seam
    $S_{K\rightarrow H'}$ is only defined  if $(H') <(K)$. Let $x\in
S_{K\rightarrow H'}\subset
    \P_{0}^{(H')}\subset \P_0$ and $U_{x}$ an open neighborhood of $x$
in $\P_{0}$.  So $\pi^{-1}(U_x)$ is an open neighborhood of a point
$z\in\J^{-1}(0)$ such that
$\pi (z)=x$, where $\pi$ denotes the orbit projection, $\pi:
\mathbf{J}^{-1}(0)\rightarrow P_0$.

 As the point $x$ projects under the map $\tau^0$ to $m\in M^{(K)}$
then we can assume  without loss of generality, that $\tau (z) = y
\in M_{(K)}$ satisfying $G_y=K$. From \eqref{pointdecomposition},
the zero momentum level set restricted to the fiber over $y$ is
given by
    $$\J_y^{-1} (0) = \J_{(K)y}^{-1}(0) \oplus N^*_yM_{(K)}.$$
Now note that because $\pi(z)=x\in S_{K\rightarrow H'}$ then
$$z\in\J_{(K)y}^{-1}(0) \times (N^*_yM_{(K)})_{(H')}$$
where the orbit type on the conormal fiber refers to the linear
$K$ action. Recall from the orbit type decomposition of the
conormal fiber $N^*_yM_{(K)}$ that for any $(H)\in I_M$ such that
$(H)< (K)$ then $(N^*_yM_{(K)})_{(H)}\neq\emptyset$, and
consequently $ (N^*_yM_{(K)})_{(H')}\subset \partial
(N^*_yM_{(K)})_{(H)}$ if $(H)<(H')<(K)$. This means that there is
a point $z'\in \pi^{-1}(U_x)\cap(\J_{(K)y}^{-1}(0) \times
(N^*_yM_{(K)})_{(H)}) $, from where {\it (v)} easily follows once
we note that $\pi(z')\in U_x\cap S_{K\rightarrow H}$.

({\it 3}) follows from the definition of a $\Sigma$-decomposed
   surjective submersion, since $\tau^0\vert_{C_L}=t^L$ and
   $\tau^0\vert_{S_{K'\rightarrow K}}=t^{K'\rightarrow K}$ are
surjective submersions  and the pieces of the coisotropic
decomposition of $\mathcal{P}_0$ are the $C_{L}$'s and the seams,
and the pieces of the orbit type decomposition of $M/G$ are
$M^{(L)}$ for every $(L)\in I_{M}$.

Finally,  ({\it 4}) is obvious from the construction of the
coisotropic and symplectic decompositions.
    \end{proof}
From the frontier conditions {\it (i)} to {\it (iii)} it is clear
that two cotangent bundles $C_K$ and $C_H$ are stitched along the
corresponding seam $S_{K\rightarrow H}$. The pieces of the
coisotropic decomposition are in one-to-one correspondence with
the connectable pairs of $I_M$, where to a connectable pair of two
copies of a same class $H\rightarrow H$ corresponds the cotangent
bundle $C_H$, and for different classes $K\rightarrow H$, $(H)\neq
(K)$ the corresponding piece is a seam $S_{K\rightarrow H}$. Thus
Theorem \ref{coisotropicdecomposition} allows us to obtain the
coisotropic decomposition lattice with only the knowledge of the
lattice $I_M$.

\section{From $\Sigma$-decompositions to stratifications}
It was the objective of this paper to give a description of the topology
and geometry of the reduced space $\P_0$, and for
a number of important reasons such a description based in the stratified nature of the
singular spaces involved is more
desirable than the one based only in the weaker concept of
$\Sigma$-decompositions. In this section we upgrade our previous
topological
results and in the following we will concentrate on giving meaning and
justification to the following assertion:
\begin{theorem}\label{theupgrade}All the $\Sigma$-decomposed spaces in
Theorems \ref{thm8} and \ref{coisotropicdecomposition} are
stratified spaces with the unique stratifications induced by their
$\Sigma$-decompositions. Consequently all the maps involved are
morphisms of stratified spaces. In particular $\tau^0$ and
$\tau^L$ are stratified surjective submersions.
\end{theorem}
We need then an appropriate definition of stratification and
morphism of stratified spaces. We will follow closely the
reference \cite{Pflaum} for the definitions in the rest of the
section. We caution the reader that other authors use different
definitions for the same terminology (for example the definition
of stratification found in \cite{SaLe}, which also includes the
extra properties of being a cone space). Let $X$ be a topological
space and $\mathcal{S}$ a map that associates to each point $x\in
X$ the set germ $\mathcal{S}_x$ at $x$ of a locally closed subset
of $X$. Recall that the set germ of a set $A$ at $x\in A$ is the
equivalence class $[A]_x$ of $A$ at $x$ defined by $[A]_x=[B]_x$
if both $A$ and $B$ are subsets of $X$ containing $x$ and such
that there exists an open neighborhood $U$ of $x$ satisfying
$A\cap U=B\cap U$.

From now on we shall call a $\Sigma$-decomposition for which,
given any piece, all its connected components have the same
dimension, a {\em decomposition}.
\begin{definition} In the previous conditions, the map $\mathcal{S}$ is
said to be a stratification of $X$ if for any point $x\in X$,
there exists an open neighborhood $U$ containing $x$ and a
decomposition $\mathcal{Z}$ of $U$ satisfying:  For any $y\in U$,
$\mathcal{S}_y=[Z]_y$, with $Z\in\mathcal{Z}$ the piece containing
$y$. The pair $(X,\mathcal{S})$ is called a stratified space.

Let $(X,\mathcal{S})$ and $(Y,\mathcal{T})$ be two stratified
spaces and $f:X\rightarrow Y$ a continuous map between the
underlying topological spaces. $f$ is called a morphism of
stratified spaces if for every $x\in X$ there exist neighborhoods
$V$ of $f(x)$ and $U\subset f^{-1}(V)$ of $x$ with decompositions
$\mathcal{X}$ and $\mathcal{Y}$ inducing $\mathcal{S}\vert_U$ and
$\mathcal{T}\vert_V$ res\-pec\-ti\-ve\-ly, such that for every
$x'\in U$ there is an open neighborhood $W\subset U$ containing
$x'$ such that the restriction $f\vert_W$ maps the intersection of
the piece $S$ containing $x'$ with $W$ into a piece
$R\in\mathcal{Y}$ and $f\vert_{S\cap W}:S\rightarrow R$ is smooth.

We will say that $f$ is a stratified immersion (resp. submersion,
diffeomorphism, etc) if so are all the restrictions $f\vert_{S\cap
W}$ at every point $x\in X$.
\end{definition}

Obviously if $(X,\mathcal{X})$ is a decomposed space, for any
neighborhood $U$ of any point, $(U,\mathcal{X}\vert_U)$ is again a
decomposed space, and then we can give $X$ the structure of a
stratified space associating to each of its points $x$ the set
germ of the piece containing $x$. This stratification is said to
be induced by the decomposition $\mathcal{X}$. As an immediate
consequence a morphism of decomposed spaces is a morphism of the
induced stratified spaces.

A $\Sigma$-decomposition $\mathcal{X}$ in principle does not
induce a stratification, since $\mathcal{X}\vert_U$ could be a
$\Sigma$-decomposition instead of a decomposition of $U$ no matter
how $U$ is chosen as we can see in the following example: Consider
the subspace $X$ of $\mathbb{R}^3$ given by the $(x_1,x_2)$-plane
and the $x_3$-axis. Let $X_1=X\backslash \mathbf{0},\,
X_2=\mathbf{0}$. Obviously $X_1$ and $X_2$ are $\Sigma$-manifolds
and the partition $X=X_1\cup X_2$ is a $\Sigma$-decomposition of
$X$, but for any open neighborhood $U$ of $\mathbf{0}$ the induced
partition of $U$ is again a $\Sigma$-decomposition, so the map
associating to each point the equivalence class of the piece
containing it is not a stratification.

However,  in the special case of the orbit type
$\Sigma$-decomposition of a proper $G$-manifold $M$ it is possible
to induce a decomposition of a suitable open neighborhood of an
arbitrary point. Furthermore, it is possible to guarantee that the
secondary and coisotropic $\Sigma$-decompositions induced from the
orbit type one are locally decompositions, fulfilling the
requirements for inducing stratifications, for which the
decomposed morphisms are automatically stratified morphisms.

The reason for this lies once again in the local model of an
invariant neighborhood $U$ of an orbit $G\cdot m$ given by the
tubular neighborhood $G\times_H S_m$ where $G_m=H$ and $S_m$ is a
linear slice orthogonal to the directions tangent to the orbit at
$m$. In this model the orbit type $U_{(L)}$ is represented by
$G\times_H (S_m)_{(L)}$, where $L$ must be a subgroup of $H$ and
the action on the linear slice is the linear $H$-action by
isometries with respect to the restriction of the inner product in
$T_mM$. But it is known that the partition of a vector space by
orbit types with respect to the linear representation of a compact
Lie group is a decomposition (see for instance Lemma 4.10.12 of
\cite{Chossat}). Consequently the induced $\Sigma$-decomposition
of $U$, consisting of the intersection of pieces in $M$ with $U$
is actually a decomposition since the pieces are of the form
$G\times_H (S_m)_{(L)}$, having its connected components the same
dimension.

To see that the coisotropic decomposition is a stratification,
first recall that the map $\tau^0:\P_0\rightarrow M/G$ is an open,
$\Sigma$-decomposed map. Now, choosing a suitable small enough
open set, $U$ in $\P_0$, it will project to a decomposed open set,
$O:=\tau^0(U)$ where all pieces have components of the same
dimension, since $M/G$ is locally decomposed. A connected
component of $S_{H\rightarrow L}\cap U$ projects under $\tau^0$ to
a connected component of $M^{(H)}\cap O$, and its dimension is
determined by the dimension of this component of $M^{(H)}\cap O$
and the dimension of some other connected component of
$M^{(L)}\cap O$ as we have seen in the proof of
Theorem~\ref{secondsymp}. Since $O$ is a decomposed space then all
these pieces of the form $M^{(L)}\cap O$ have the same dimension,
from where it follows that all the connected components of
$S_{H\rightarrow L}\cap U$ have the same dimension, and therefore
$U$ is a decomposed open set in $\P_0$, proving that the
coisotropic decomposition is a stratification. Similar arguments
work for the secondary decomposition, and so we conclude
Theorem~\ref{theupgrade}. We are therefore justified to use the
terminology secondary and coisotropic stratifications, as well as
their corresponding stratification lattices.

\section{An example}
We will illustrate the main results obtained in  this paper with
an example that is simple, yet rich enough to
show the extra structure appearing in singular symplectic reduction for
cotangent bundles. We will compute the
secondary and coisotropic stratifications exhibiting explicitly the
corresponding frontier conditions predicted in Theorems \ref{thm8} and
\ref{coisotropicdecomposition}.

Consider the $G=\mathbb{Z}_2\times S^1$ action on
$M=\mathbb{R}^3$, where  $S^1$ acts by rotations around the
$x_3$-axis and $\mathbb{Z}_2$ by reflections with respect to the
plane $(x_1,x_2)$. The isotropy lattice and the decomposition
lattice for this action are shown in Figure~\ref{examples}. Let
$\mathbb{R}^3$ be equipped with the Euclidean inner product  which
defines a $G$-invariant Riemannian metric for this action.
Identifying $T^*\mathbb{R}^3$ with $\mathbb{R}^3\times
\mathbb{R}^3$ then the cotangent lifted action is   diagonal,
$g\cdot(v_1,v_2)=(g\cdot v_1,g\cdot v_2)$ for $g\in G$ and
$v_1,v_2\in \Rb^3$. Let $(x_1,x_2,x_3,y_1,y_2,y_3)$ be the
coordinates of  the vector $(\mathbf{x},\mathbf{y})\in
\mathbb{R}^3\times \mathbb{R}^3$ with respect to the canonical
basis.

The  ring of $G$-invariant polynomials,
$\P^G(\mathbb{R}^3\times \mathbb{R}^3)$, is generated by
$$\begin{array}{lclclcl}
\sigma_1 & = & x_1^2+x_2^2 + y_1^2+ y_2^2,&\quad & \rho_1 & = &
x_3^2+y_3^2, \\
\sigma_2 & = & 2(x_1y_1+x_2y_2), &\quad & \rho_2 & = & 2x_3y_3,  \\
\sigma_3 & = & y_1^2+y_2^2-x_1^2-x_2^2, &\quad & \rho_3 & = &
y_3^2-x_3^2,  \\
j & = & x_1y_2-x_2y_1. & & & &
\end{array}$$
These polynomials are subject to the relations
$$
\sigma_1  \geq  0,\quad \rho_1  \geq  0, \quad
\sigma_1^2  =  \sigma_2^2+\sigma_3^2 + 4 j^2, \quad \rho_1^2  =
\rho_2^2+\rho_3^2.  \\
$$
Note
  that the relations between the $\sigma$'s and the $\rho$'s are
uncoupled if   $j$ is zero. The   momentum map for the cotangent
lifted action of $G$ is  $\J (\mathbf{x},\mathbf{y}) = j$.

Let now $Z$ be a $G$-invariant subset of
$\mathbb{R}^3\times\mathbb{R}^3$ such
that $j$ is constant on $Z$. Consider two copies of
$\mathbb{R}^3$, which will be denoted by $\mathbb{R}^{3\sigma}$ and
$\mathbb{R}^{3\rho}$ and the maps
$\chi_\sigma:Z\rightarrow\mathbb{R}^{3\sigma}$ and
$\chi_\rho:Z\rightarrow\mathbb{R}^{3\rho}$ defined as
\begin{eqnarray*}
\chi_\sigma (z)  =  (\sigma_1(z),\sigma_2(z),\sigma_3(z)), \quad
\chi_\rho (z)  =  (\rho_1(z),\rho_2(z),\rho_3(z))
\end{eqnarray*}
for every $z\in Z$. The Hilbert map
$\chi:=(\chi_\sigma,\chi_\rho):Z\rightarrow
\mathbb{R}^{3\sigma}\times \mathbb{R}^{3\rho}$ is $G$-invariant,
and due to the relation between the polynomials its image
$\mathrm{Im}\,\chi=\mathrm{Im}\,\chi_\sigma\times\mathrm{Im}\,\chi_\rho\
\in\mathbb{R}^{3\sigma}\times \mathbb{R}^{3\rho}$ is a topological
space equipped with the relative topology which is a
semi-algebraic variety.  The Tarski-Seidenberg Theorem (see
\cite{DuiKol} and references therein for a more detailed
explanation) gives that $\mathrm{Im}\,\chi$ has a canonical
(Whitney) stratification.  By invariant theory the map $\chi$
restricts to a homeomorphism $\overline{\chi}:Z/
G\rightarrow\mathrm{Im}\,\chi_\sigma\times\mathrm{Im}\,\chi_\rho\in\mathbb{R}^{3\sigma}\times
\mathbb{R}^{3\rho}$ that happens to be an isomorphism of
stratified spaces if $Z/G$ is endowed with the orbit type
stratification.  In order to apply the results obtained in
previous sections we will study  the case $Z=\J^{-1}(0)$   through
the image of $\chi$.

\begin{figure}[htbp]
\begin{center}
\includegraphics{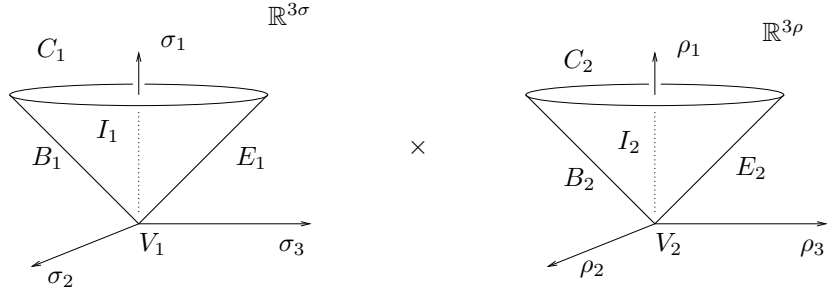}
\end{center}
\caption{\label{doublecone}The reduced singular space $\P_0$ as a
product of two cones}
          \end{figure}

The zero level set of the
momentum map is
$Z=\mathbf{J}^{-1}(0)=\{\mathbf{(x,y)}\in\mathbb{R}^6\,\vert\,
j(\mathbf{x}, \mathbf{y})=0\}$.  So we can identify $\P_0$ with the
direct
product of the two cones defined by the relations
$$C_1:\,\sigma^2_1=\sigma_2^2+\sigma_3^2,\quad\mathrm{and}\quad
C_2:\,\rho^2_1=\rho_2^2+\rho_3^2.$$ This realization of $\P_0$ is
shown in Figure~\ref{doublecone}. For future reference  in
Figure~\ref{doublecone} we mark some  subsets on each of the
cones. For instance in $C_1$ the vertex is marked as $V_1$, the
straight line $\sigma_1=\sigma_3$ excluding the origin is labelled
$E_1$, the opposite line $\sigma_1=-\sigma_3$ also except the
origin is labelled as $B_1$, and finally all the cone except
$V_1\cup E_1$ is called $I_1$ ($I_1$ contains $B_1$). Note from
the defining equation of $C_1$ that $B_1$ and $E_1$ form an angle
of $\pi/2$. Analogous definitions  apply to $C_2$.

By Proposition \ref{prop10} we know that the orbit types present
in $Z$ are exactly those which are present in $M$, i.e. the
elements of $I_M$. This implies that the symplectic strata of
$\P_0$ are in one-to-one correspondence with the strata of $M$,
and that both spaces exhibit an identical stratification lattice
as we will verify now. Indeed, studying the diagonal action
restricted to $Z$ one finds easily  the following orbit types:
$$\begin{array}{lll}
Z_{(\mathbb{Z}_2\times S^1)} & = & \left\{(\mathbf{0},\mathbf{0})\right\}\\
Z_{(\mathbb{Z}_2)} & = &
\left\{(\mathbf{x},\mathbf{y})\in\mathbb{R}^6\,\vert\
x_3=y_3=0,\, x_1y_2-x_2y_1=0)\right\},\\
Z_{(S^1)} & = &
\left\{(\mathbf{x},\mathbf{y})\in\mathbb{R}^6\,\vert\,x_1=x_2=y_1=y_2=0,\,\,(
x_3,y_3)\neq(0,0)\right\}\\
Z_{(1)} & = & Z\setminus (Z_{(\mathbb{Z}_2\times S^1)} \cup
Z_{(\mathbb{Z}_2)} \cup Z_{(S^1)})
\end{array}
$$
Using the image of the map $\chi$ we have
$$
\begin{array}{lllll}
\P_0^{(\mathbb{Z}_2\times S^1)} & = & V_1\times V_2&&\\
\P_0^{(\mathbb{Z}_2)} & = & (I_1\cup E_1) \times V_2&=&(C_1\setminus
V_1)\times V_2\\
\P_0^{(S^1)} & = & V_1 \times (I_2\cup E_2)&=&V_1\times (C_2\setminus
V_2)\\
\P_0^{(1)} & = & (I_1\cup E_1)\times (I_2\cup E_2)&=&(C_1\setminus
V_1)\times (C_2\setminus V_2).
\end{array}
$$
The above sets are the strata of the symplectic stratification
lattice predicted by Theorem~\ref{SL}.  This lattice is shown in
Figure~\ref{examplesecondary} a).

Recall that the strata for the secondary stratification of each
symplectic stratum $\P_0^{(L)}$ are of two types, cotangent
bundles $C_L$ and seams $S_{H\rightarrow L}$ with $(H) > (L)$
defined by \eqref{seam}. Let us now study  the secondary
stratification of each symplectic stratum in $\P_0$. We embed
$M=\mathbb{R}^3$ in $T^*M$ by the injection
\begin{equation}\label{exampleinclusion}(x_1,x_2,x_3)\mapsto
(x_1,x_2,x_3,0,0,0).\end{equation} We then have
$$\begin{array}{lll}
T^*M_{(\mathbb{Z}_2\times S^1)} & = &\left\{ (\mathbf{0},\mathbf{0})\right\}\\
T^*M_{(\mathbb{Z}_2)} & = & \left\{(x_1,x_2,0,y_1,y_2,0),\quad
(x_1,x_2)\neq
(0,0)\right\}\\
T^*M_{(S^1)} & = & \left\{(0,0,x_3,0,0,y_3), \quad x_3\neq 0\right\}\\
T^*M_{(\mathbf{1})} & = & \left\{(\mathbf{x},\mathbf{y}), \quad
\mathbf{x}\in
M_{(\mathbf{1})}\right\}\\
N^*M_{(\mathbb{Z}_2\times S^1)} & = &
\left\{(\mathbf{0},\mathbf{y}), \quad
\mathbf{y}\in\Rb^3\right\}\\
N^*M_{(\mathbb{Z}_2)}  & = & \left\{(x_1,x_2,0,0,0,y_3), \quad
(x_1,x_2)\neq
(0,0)\right\}\\
N^*M_{(S^1)} & = & \left\{(0,0,x_3,y_1,y_2,0), \quad x_3\neq 0\right\}\\
N^*M_{(\mathbf{1})} & = & \left\{(\mathbf{x},\mathbf{0}), \quad
\mathbf{x}\in M_{(\mathbf{1})}\right\}.
\end{array}
$$
Computing the seams and the cotangent bundles we obtain the
following re\-a\-li\-za\-tion of these two types of pieces in the
image of $\chi$:
$$\begin{array}{lllllll}
C_{\mathbb{Z}_2\times S^1} & = & V_1\times V_2 & \quad &
S_{\mathbb{Z}_2\times S^1\rightarrow \mathbb{Z}_2} & = & E_1\times V_2\\
C_{\mathbb{Z}_2} & = & I_1 \times V_2 & \quad & S_{\mathbb{Z}_2\times
S^1\rightarrow S^1} & = & V_1\times E_2\\
C_{S^1} & = & V_1\times I_2 & \quad & S_{\mathbb{Z}_2\times
S^1\rightarrow
\mathbf{1}} & = & E_1\times E_2\\
C_{\mathbf{1}} & = & I_1\times I_2 & \quad & S_{\mathbb{Z}_2\rightarrow
\mathbf{1}} & =  & I_1\times E_2\\
         &   &               &       & S_{S^1\rightarrow\mathbf{1}} & = &
E_1\times I_2.
\end{array}$$
According to the results of Theorem \ref{thm8}
  the secondary stratifications of the symplectic strata are:
$$\begin{array}{lll}
\P_0^{(\mathbb{Z}_2\times S^1)} & = & C_{(\mathbb{Z}_2\times S^1)}
=V_1\times V_2\\
\P_0^{(\mathbb{Z}_2)} & = & C_{\mathbb{Z}_2} \cup
S_{\mathbb{Z}_2\times S^1\rightarrow \mathbb{Z}_2}  =  (I_1\times
V_2)\cup  (E_1\times
V_2)\\
\P_0^{(S^1)} & = & C_{S^1}\cup S_{\mathbb{Z}_2\times
S^1\rightarrow S^1} =
 (V_1\times I_2)\cup (V_1\times E_2)\\
\P_0^{(\mathbf{1})} & = & C_{\mathbf{1}}\cup
S_{\mathbb{Z}_2\rightarrow \mathbf{1}}\cup
S_{S^1\rightarrow\mathbf{1}}\cup S_{\mathbb{Z}_2\times
S^1\rightarrow \mathbf{1}}  =  (I_1\times I_2)\cup
(I_1\times E_2)\\
& & \hspace{4.5cm} \cup  (E_1\times I_2)\cup (E_1\times E_2).
\end{array}$$
The corresponding stratification lattices are shown in
Figure~\ref{examplesecondary} (b)-(e). The coisotropic
stratification lattice is shown in
Figure~\ref{examplecoisotropic}. These lattices are constructed
using the
  results of  Theorems~\ref{thm8} and \ref{coisotropicdecomposition},
  and the corresponding frontier conditions can be verified from the above
  expressions.
        \begin{figure}[htbp]
   $$\begin{array}{cc}
\xymatrix{ & \P_0^{(\mathbb{Z}_2\times S^1)}\ar[dl]\ar[dr] &\\
\P_0^{(\mathbb{Z}_2)}\ar[dr] & & \P_0^{(S^1)}\ar[dl]\\
  & \P_0^{(\mathbf{1})} &
  \\}
&
\xymatrix{ & S_{\mathbb{Z}_2\times
S^1\rightarrow \mathbb{Z}_2}\ar[d]\\
C_{\mathbb{Z}_2\times S^1}  &  C_{\mathbb{Z}_2}\\
  }\\
&\\
a)&b)\qquad\qquad\qquad c)\\
&\\
\xymatrix{S_{\mathbb{Z}_2\times S^1\rightarrow S^1}\ar[d]\\
C_{S^1}\\} &
\xymatrix{ &  S_{\mathbb{Z}_2\times S^1\rightarrow 1}\ar[dl]\ar[dr] &\\
S_{\mathbb{Z}_2\rightarrow \mathbf{1}}\ar[dr] & &
S_{S^1\rightarrow\mathbf{1}}\ar[dl]\\
  & C_{\mathbf{1}} &\\}\\
&\\
d)&e)\\
\end{array}$$\caption{\label{examplesecondary} a) Symplectic
stratification of $\P_0$. Secondary stratifications of: b)
$\P_0^{(\mathbb{Z}_2\times S^1)}$,
c) $\P_0^{(\mathbb{Z}_2)}$, d) $\P_0^{S^1}$ and e)
$\P_0^{(\mathbf{1})}$.}
              \end{figure}
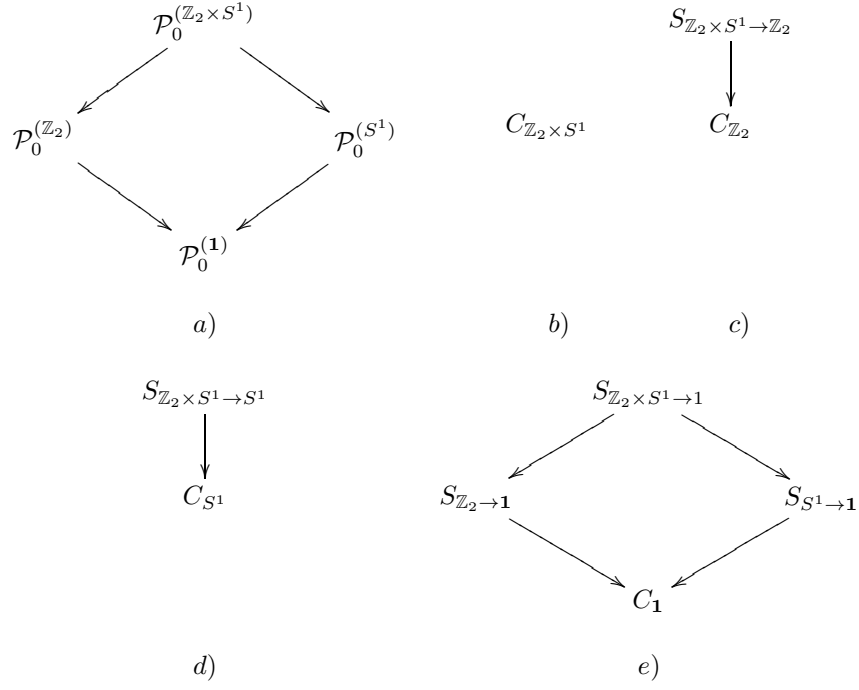
              \begin{figure}[htbp]
              $$\xymatrix{ & C_{\mathbb{Z}_2\times S^1}\ar[dl]\ar[dr] &
\\
S_{\mathbb{Z}_2\times S^1\rightarrow \mathbb{Z}_2}\ar[d]\ar[dr] & &
S_{\mathbb{Z}_2\times S^1\rightarrow S^1}\ar[d]\ar[dl]\\
C_{\mathbb{Z}_2}\ar[d] & S_{\mathbb{Z}_2\times S^1\rightarrow
\mathbf{1}}\ar[dl]\ar[dr]  & C_{S^1}\ar[d]\\
S_{\mathbb{Z}_2\rightarrow \mathbf{1}}\ar[dr] & &
S_{S^1\rightarrow\mathbf{1}}\ar[dl]\\
& C_{\mathbf{1}} &
}$$
\caption{\label{examplecoisotropic} Coisotropic stratification of
$\P_0$.}
         \end{figure}
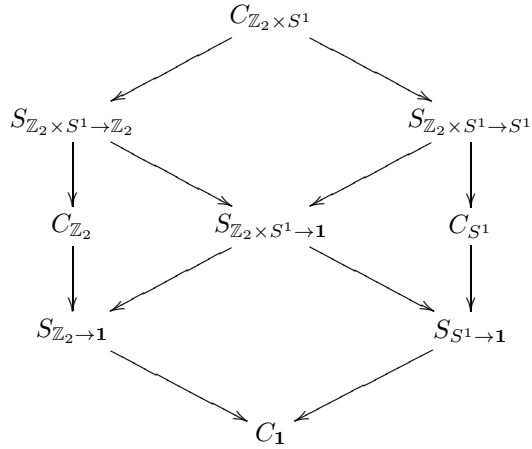
We describe now the bundle structure of these stratifications:
using equation~\eqref{exampleinclusion} we realize the quotient
$M/G$ as the subset of the image of $\chi$ given by $(B_1\cup
V_1)\times (B_2\cup V_2)$. The corresponding strata of its orbit
type stratification are:
$$\begin{array}{ccccccc} M^{(\mathbb{Z}_2\times S^1)} & = &
V_1\times V_2 & & M^{(\mathbb{Z}_2)} &
= & B_1\times V_2\\
M^{(S^1)} & = & V_1 \times B_2 & & M^{(\mathbf{1})} & = &
B_1\times B_2\\
\end{array}$$
The map $\tau^0:\P_0\rightarrow M/G$ is obtained as follows. Let
$z=(\mathbf{x}_1,\mathbf{x}_2)\in C_1\times C_2$ be a point of
$\P_0$, then $\tau^0(z)$ is a point
$(\mathbf{b}_1,\mathbf{b}_2)\in B_1\times B_2$ where
$\mathbf{b}_1$ is the point in the intersection of $B_1$ and the
unique parabola obtained by sectioning the cone $C_1$ with a plane
orthogonal to $B_1$ at $\mathbf{x}_1$. Analogously one defines in
this way the point $\mathbf{b}_2\in C_2$.

\section{Final Remarks}
We have studied the global picture of two new stratifications of
the zero momentum singular reduced space for a cotangent lifted
action. The results obtained raise several natural questions which
have not been addressed in this work.

First, as mentioned in the Introduction, it would be interesting
to determine if these reduced spaces, together with the secondary
and coisotropic stratifications,  have  conical structure, satisfy
Whitney conditions and/or admit singular atlases and smooth
structures, as it happens for the symplectic stratification (see
\cite{SaLe} and \cite{OrRa2003}). For this, the Symplectic Slice
Theorem of Marle, Guillemin and Sternberg is too weak, and a
cotangent bundle adapted version of it, which could detect the
secondary and coisotropic strata would be needed. Unfortunately,
such a technology does not yet exist in full generality. Some
steps have been done in \cite{Schmah} but it is still lacking a
general result. We expect however that advances in this line of
research can lead to the proof that the secondary and the
coisotropic stratifications enjoy the conical property and satisfy
Whitney conditions like the symplectic stratification.

A different direction of study consists of describing reduction at
nonzero momentum. At least for reduction at momentum values with
trivial coadjoint orbits it is also possible to obtain a secondary
and coisotropic stratification with some modifications of the
technology used here. This will appear elsewhere. For general
momenta the problem is much more involved since the coadjoint
representation interacts with the action on the base manifold to
produce an isotropy lattice of the momentum level set
$\mathbf{J}^{-1}(\mu)$. These are aspects of ongoing work on the
subject.

Even when the secondary stratification of each symplectic stratum
is not invariant for the reduced Hamiltonian flow, the fact that
it captures the bundle structure of the reduced space might be
useful for understanding certain qualitative aspects of the
reduced dynamics. A typical situation would be described by the
following observation: Even when it is known that the isotropy is
preserved by the reduced dynamics (and hence the symplectic strata
are dynamically invariant), this is not true for the isotropy of
the projected dynamics onto $M/G$. It is precisely when the
Hamiltonian evolution crosses different seams within its ambient
symplectic stratum that changes in the isotropies of the base
points occur. Natural questions arise then, like under what
conditions symmetric Hamiltonian flows preserve isotropy both in
phase and configuration spaces, or when a Hamiltonian evolution is
maximal, in the sense that every secondary stratum (and hence
every possible isotropy type in the base) is crossed. This will be
the object of further research.
\section*{Acknowledgements}
This work was partially support by  the  EU funding for the
Research Training Network MASIE,  Contract No. HPRN-CT-2000-00113
and by FCT (Portugal) through the programs POCTI/FEDER. We would
like to thank Mark Roberts for pointing out a mistake in an early
stage of this work and Tanya Schmah for several useful suggestions
and comments.

\end{document}